\documentclass[10pt]{article}
\usepackage[left=1in,top=1.25in,right=1in,bottom=1.25in]{geometry} %adjust margins here
\usepackage{amssymb,amsmath,amsthm}
\numberwithin{equation}{section}
\usepackage{color}
\usepackage{bbm} %for the indicator function symbol
\usepackage{accents} %for the vector symbol
\makeatletter
\newcommand{\subjclass}[2][2020]{%
	\let\@oldtitle\@title%
	\gdef\@title{\@oldtitle\footnotetext{#1 \emph{Mathematics subject classification.} #2}}%
}
\newcommand{\keywords}[1]{%
	\let\@@oldtitle\@title%
	\gdef\@title{\@@oldtitle\footnotetext{\emph{Key words and phrases.} #1.}}%
}
\makeatother
\newtheorem{thm}{Theorem}[section]
\newtheorem{cor}[thm]{Corollary}
\newtheorem{prop}[thm]{Proposition}
\newtheorem{lem}[thm]{Lemma}
\theoremstyle{definition}
\newtheorem{defn}[thm]{Def{i}nition}
\newtheorem{rmk}[thm]{Remark}

\newcommand{\pf}{\noindent \textit{Proof}:\ }
\newcommand{\R}{\mathbb{R}}
\newcommand{\C}{\mathbb{C}}
\newcommand{\N}{\mathbb{N}}
\newcommand{\E}{\mathbb{E}}
\newcommand{\de}{\,\text{\rm d}}
\renewcommand{\qed}{\hfill{\ensuremath\square}}

\DeclareMathOperator*{\prlim}{prlim}
\DeclareMathOperator*{\indlim}{indlim}

\makeatletter %for the double left/right angle symbol
\newsavebox{\@brx}
\newcommand{\llangle}[1][]{\savebox{\@brx}{\(\m@th{#1\langle}\)}%
	\mathopen{\copy\@brx\kern-0.5\wd\@brx\usebox{\@brx}}}
\newcommand{\rrangle}[1][]{\savebox{\@brx}{\(\m@th{#1\rangle}\)}%
	\mathclose{\copy\@brx\kern-0.5\wd\@brx\usebox{\@brx}}}
\makeatother
\usepackage{hyperref} %for hyperlinks; this should be last in the preamble 

\begin{document}
	
\title{Stochastic analysis for vector-valued generalized grey Brownian motion}
\author{Wolfgang Bock, Martin Grothaus, Karlo Orge}
\date{ }
\maketitle

\begin{abstract}
	In this article, we show that the standard vector-valued generalization of a generalized grey Brownian motion (ggBm) has independent components if and only if it is a fractional Brownian motion. In order to extend ggBm with independent components, we introduce a vector-valued generalized grey Brownian motion (vggBm). The characteristic function of the corresponding measure is introduced as the product of the characteristic functions of the one-dimensional case. We show that for this measure, the Appell system and a calculus of generalized functions or distributions are accessible. We characterize these distributions with suitable transformations and give a $d$-dimensional Donsker's delta function as an example for such distributions. From there, we show the existence of local times and self-intersection local times of vggBm as distributions under some constraints, and compute their corresponding generalized expectations. At the end, we solve a system of linear SDEs driven by a vggBm noise in $d$ dimensions. 
\end{abstract}

	\section{Introduction}
	\label{sec:introduction}
Many applications call for processes with long-range dependence and complex correlation structures. As a generalization of Brownian motion, fractional Brownian motion (fBm) is used to model such dynamics, based on its correlated increments, which imply short or and long-range dependence \cite{Mis08,biagini2008,nourdin}.
Fractional Brownian motion is neither a semi-martingale nor a Markov process, except for the Brownian motion case. Hence, it is not accessible by standard stochastic calculus, and thus challenging from the mathematical point of view.  
There are various ways to cast fBm into the classical Brownian motion framework, starting with the famous definition by Mandelbrot and van Ness \cite{MandelbrotNess1968}. This idea is also the starting point for a characterization of fBm using an infinite superposition of Ornstein-Uhlenbeck processes with respect to the standard Wiener process; compare the works of Carmona, Coutin, Montseny, and Muravlev \cite{CCM2000,CCM2003,MU2011} or also the monograph of \cite{Mis08}. Recently, further applications of this representation have for instance been investigated in \cite{HS2019} with a focus on finance and in \cite{BD2020} in the context of optimal portfolios.

The Mandelbrot-van Ness representation can be used to represent fBm in the framework of white noise analysis \cite{Mis08,biagini2008,nourdin, Ben03}.
White noise analysis has evolved into an infinite dimensional distribution theory, with rapid developments in mathematical structure and applications in various domains; see, e.g.~the monographs \cite{HKPS93, Oba94, Kuo96}. Various characterization theorems \cite{PS91, KLPSW96, GKS97, HS17, GMN21} are proven to build up a strong analytical foundation.
Almost at the same time, first attempts were made to introduce a non-Gaussian infinite dimensional analysis, by transferring properties of the Gaussian measure to the Poisson measure \cite{I88} with the help of bi-orthogonal generalized Appell systems \cite{Da91, ADKS96, KSWY98}. This approach is suitable for many measures, like the Gaussian measure and the Poisson measure \cite{KSS97}. 
Mittag-Leffler Analysis is established in \cite{GJRS15} and \cite{GJ16}. In fact, it generalizes methods from white noise calculus to the case, where in the characteristic function of the Gaussian measure the exponential function is replaced by a Mittag-Leffler function. The corresponding grey Brownian motion is in general neither a martingale nor a Markov process. Moreover, it is not a possible - as in the Gaussian case - to find a proper orthonormal system of polynomials for the test and generalized functions. Here, it is necessary to make use of the aforementioned Appell system of bi-orthogonal polynomials.
The grey noise measure \cite{Sch92, MM09} is included as a special case in the class of Mittag-Leffler measures, which offers the possibility to apply the Mittag-Leffler analysis to fractional differential equations, in particular to fractional diffusion equations \cite{Sch90, Sch92}, which carry numerous applications in science, like relaxation type differential equations or viscoelasticity. 
Fractional processes were motivated by phenomena in heterogeneous media modeled by fractional partial differential equations; see \cite{kochubei2004,mainardi2010,meerschaert2004limit}. Corresponding stochastic processes governed by these equations have applications in science, engineering and finance \cite{rangarajan2003processes,magdziarz2009black,metzler2000random,metzler2004restaurant,scalas2006five}. 
Fractional time derivatives are used to model sticking of particles in porous media \cite{schumer2003fractal}. In statistical physics, fractional time derivatives reflects random waiting times between particle jumps \cite{meerschaert2004limit}.  Detailed discussion of such processes is also found in \cite{bender2021stochastic}.  For a detailed study of the special class of heavy tailed processes, see \cite{leonenko2017heavy,leonenko2020approximation,heyde2005student}. An approach using subordination can be found in \cite{bender2021stochastic}.

With the help of Mittag-Leffler analysis, a relation between the fractional heat equation  and the associated process, i.e.~grey Brownian motion, was proven in \cite{GJ16}.  In \cite{BS17}, Wick-type stochastic differential equations and Ornstein-Uhlenbeck processes were solved within the framework of Mittag-Leffler analysis. In \cite{bock2020integral}, the results of of \cite{MU2011} and \cite{HS2019} for fBm were extended to the non-Gaussian case of ggBm by representing it via generalized grey Ornstein-Uhlenbeck processes, using that ggBm can be written as a product of a positive and time-independent random variable and a fBm \cite{mainardi2010functions}. A similar representation can be found in \cite{d2018centre}.
While recent progress is visible for the one-dimensional generalized grey Brownian motion, many results can not be carried over directly to a multi-dimensional case. Moreover, properties as independent components in a multi-dimensional process are desireable in many applications. 	

In this article, we study the generalization of the generalized grey Brownian motion with componentwise independence. This process is not the $d$-dimensional generalized grey Brownian motion for $d\geq 2$, since we show that the componentwise independence holds only in the Gaussian case. For this process, we study the accessibility to Appell systems. We establish a characterization and tools for the analysis of the corresponding distribution spaces. Moreover, we give explicit examples of the use of this characterization and the tools on Donsker's delta function, local times, self-intersection local times, and $d$-dimensional stochastic differential equations driven by a vggBm noise.
	
	\section{Preliminaries}
	This section provides an overview of the Mittag-Leffler measure defined on the dual space of a nuclear Frech\'{e}t space from \cite{GJRS15}.
	\subsection{Nuclear triples}
	Let $\mathcal{H}$ be a real separable Hilbert space with scalar product and induced norm denoted by $(\cdot,\cdot)$ and $|\cdot|$, respectively. Let $\mathcal{N}$ be a nuclear Frech\'{e}t space topologically and densely embedded in $\mathcal{H}$. Together with its dual space $\mathcal{N}'$, we obtain the following inclusions, called a nuclear (Gel'fand) triple:
	\begin{align*}
		\mathcal{N} \subset \mathcal{H} \subset \mathcal{N}',
	\end{align*}
	where we identify $\mathcal{H}$ with $\mathcal{H}'$ via the Riesz isomorphism. Without loss of generality, we assume that the nuclear Frech\'{e}t space $\mathcal{N}$ can be represented by a countable family of real separable Hilbert spaces $(\mathcal{H}_p)_{p\in\N}$ with the following properties.
	\begin{enumerate}
		\item[(N1)] For all $p\in\N_0$, the inclusion $\mathcal{H}_{p+1} \subset \mathcal{H}_{p}$ is a Hilbert-Schmidt operator.
		\item[(N2)] For all $p\in\N_0$, $|\cdot|_{p}\leq |\cdot|_{p+1}$ on $\mathcal{H}_{p+1}$.
	\end{enumerate}
	(Here, $\N_0:=\N\cup\{0\}$, $\mathcal{H}_0:=\mathcal{H}$, and for $p\in\N_0$, $|\cdot|_p$ is the induced norm of $\mathcal{H}_p$.) The space $\mathcal{N}$ is then assumed to be the projective limit of the spaces $(\mathcal{H}_p)_{p\in\N}$. That is, as a set $\mathcal{N} = \bigcap_{p\in\N} \mathcal{H}_p$, endowed with the coarsest topology such that the embeddings $\mathcal{N}\subset\mathcal{H}_p$, $p\in\N$ are continuous. Let $\mathcal{H}_{-p}$ be the dual space of $\mathcal{H}_{p}$, with corresponding norm $|\cdot|_{-p}$. By general duality theory (see, e.g.~\cite{GV64}), the dual space $\mathcal{N}'$ is equal to the inductive limit of $(\mathcal{H}_{-p})_{p\in\N}$. That is, $\mathcal{N}' = \bigcup_{p\in\N} \mathcal{H}_{-p}$, equipped with the finest locally convex topology such that $\mathcal{H}_{-p}$ is continuously embedded in $\mathcal{N}'$. It is known that this topology on $\mathcal{N}'$ is reflexive, so that it coincides with the strong topology on $\mathcal{N}'$ (see, e.g.~Appendix A.5 in \cite{HKPS93}).
	
	The canonical dual pairing between $\mathcal{N}'$ and $\mathcal{N}$ is denoted by $\langle\cdot,\cdot\rangle$, and is given as an extension of the scalar product on $\mathcal{H}$ by
	\begin{align*}
		\langle f,\varphi\rangle = (f,\varphi),\qquad f\in \mathcal{H},\ \varphi\in\mathcal{N}.
	\end{align*}
	The same notation is used for the dual pairing between $\mathcal{H}_{-p}$ and $\mathcal{H}_{p}$.
	
	The Hilbert tensor powers of $\mathcal{H}_{p}$ is denoted by $\mathcal{H}_{p}^{\otimes n}$, while the subspace of symmetric elements of $\mathcal{H}^{\otimes n}$ is denoted by $\mathcal{H}^{\widehat{\otimes} n}$. The same notations are used for $\mathcal{H}_{-p}$. The norms on $\mathcal{H}_{p}^{\otimes n}$ and $\mathcal{H}_{-p}^{\otimes n}$ are also denoted by $|\cdot|_{p}$ and $|\cdot|_{-p}$, respectively. Then $\mathcal{H}_{-p}^{\otimes n}$ is the dual space of $\mathcal{H}_{p}^{\otimes n}$ with respect to $\mathcal{H}^{\otimes n}$. The tensor powers $\mathcal{N}^{\otimes n}$ of $\mathcal{N}$ are defined as the projective limit of the spaces $(\mathcal{H}_{p}^{\otimes n})_{p\in\N}$, so that $(\mathcal{N}^{\otimes n})'$ is the inductive limit of $(\mathcal{H}_{-p}^{\otimes n})_{p\in\N}$. The symmetric tensor powers $\mathcal{N}^{\widehat{\otimes} n}$ are also defined similarly. The same notation $\langle\cdot,\cdot\rangle$ is used for all dual pairings for these tensor powers.
	
	We also use the complexification of all the real spaces described previously, denoted by a subscript $\C$. The element $f:=[f_1,f_2]$ in a complexification is denoted by $f=f_1 + if_2$. In the case of the complexification of the dual space, the dual pairing is extended in a bilinear way:
	\begin{align*}
		\langle F_1 + iF_2 , f_1 + if_2\rangle := \langle F_1,f_1\rangle - \langle F_2,f_2\rangle + i(\langle F_1,f_2\rangle + \langle F_2,f_1\rangle).
	\end{align*}
	Any linear map $L$ between two real linear spaces extends to a linear map between their corresponding complexifications, denoted by the same symbol, in a natural way:
	\begin{align*}
		L(f_1 + if_2) := L(f_1) + iL(f_2).
	\end{align*}
	The same symbols are used for the scalar product and induced norm of the complexification of a real Hilbert space, defined by
	\begin{align*}
		(f_1+if_2,g_1+ig_2) &:= (f_1,g_1) + (f_2,g_2) + i\big((f_2,g_1)-(f_1,g_2)\big),\\
		|f_1+if_2| &:= (f_1+if_2,f_1+if_2)^{1/2}.
	\end{align*}
	
	An example of a nuclear triple used in Gaussian analysis is the following: let $\mathcal{H}:=L^2(\R)$, the real Hilbert space of square-integrable functions on $\R$ with respect to the Lebesgue measure, and $\mathcal{N}:=\mathcal{S}(\R)$, the nuclear space of Schwartz test functions on $\R$. It is known that $\mathcal{S}(\R)$ is dense in $L^2(\R)$. Furthermore, $\mathcal{S}(\R)$ is the projective limit of a sequence $(\mathcal{H}_p)_{p\in\N}$ of real separable Hilbert spaces satisfying (N1) and (N2) (see, e.g., Appendix A.5 in~\cite{HKPS93}, for the explicit construction of $\mathcal{H}_p$), so that $\mathcal{S}(\R)$ is a nuclear Frech\'{e}t space. Together with the dual space $\mathcal{N}' = \mathcal{S}'(\R)$ of tempered distributions on $\R$, we obtain the following nuclear triple:
	\begin{align*}
		\mathcal{S}(\R) \subset L^2(\R) \subset \mathcal{S}'(\R).
	\end{align*}
	
	\subsection{The Mittag-Leffler measure}
	Given a nuclear triple $\mathcal{N}\subset\mathcal{H}\subset\mathcal{N}'$, the space $\mathcal{N}'$ is equipped with the $\sigma$-algebra $\mathcal{C}_\sigma(\mathcal{N}')$ generated by the cylinder sets
	\begin{align*}
		\{\omega\in\mathcal{N}' : (\langle\omega,\varphi_1\rangle, \dots, \langle\omega,\varphi_n\rangle)\in A \},
	\end{align*}
	where $n\in\N$, $\varphi_1,\dots,\varphi_n\in\mathcal{N}$ and $A\in\mathcal{B}(\R^n)$, the Borel $\sigma$-algebra over $\R^n$. Since $\mathcal{N}$ is a projective limit of a countable number of Hilbert spaces, $\mathcal{C}_\sigma(\mathcal{N}')$ coincides with the Borel $\sigma$-algebra generated by the weak and strong topologies on $\mathcal{N}'$ (see~\cite{BK95}).
	
	The definition of the Mittag-Leffler measure on $(\mathcal{N}',\mathcal{C}_\sigma(\mathcal{N}'))$ relies on the following function introduced by Mittag-Leffler in a series of papers \cite{ML03,ML04,ML05}; see also \cite{Wim05a,Wim05b}. We also introduce its generalization first appeared in~\cite{Wim05a}.
	\begin{defn}
		For $0<\beta<\infty$, the \textit{Mittag-Leffler function} $E_\beta$ is an entire function defined by its power series
		\begin{align}\label{e:mlf}
			E_\beta(z):= \sum_{n=0}^{\infty} \frac{z^n}{\Gamma(\beta n + 1)},\qquad z\in\C,
		\end{align}
		where $\Gamma$ is the Gamma function.
		In addition, for $0<\rho<\infty$, the \textit{generalized Mittag-Leffler function} $E_{\beta,\rho}$ is an entire function defined by the power series
		\begin{align*}
			E_{\beta,\rho}(z):= \sum_{n=0}^{\infty} \frac{z^n}{\Gamma(\beta n + \rho)},\qquad z\in\C.
		\end{align*}
	\end{defn}
	Note that for any $z\in\C$, $E_{\beta,1}(z) = E_\beta(z)$ and $E_{1}(z) = e^z$. Moreover, since $E_\beta$ is entire, we can calculate its derivative by differentiating term-by-term the series in (\ref{e:mlf}), and obtain
	\begin{align}\label{e:MLder}
		\frac{\text{\rm d}}{\text{\rm d}z}E_\beta(z) = \frac{1}{\beta}E_{\beta,\beta}(z).
	\end{align}
	Also, for $0<\beta\leq 1$, the map $[0,\infty)\ni x\to E_\beta(-x)\in\R$ is completely monotonic, that is, $(-1)^n E^{(n)}_\beta(-x)\geq 0$ for all $x\geq 0$. Using this fact and \cite{Pol48}, one can show in a similar manner as that of  \cite{Sch92} that the map
	\begin{align*}
		\mathcal{N}\ni \varphi \mapsto E_\beta\left(-\frac{1}{2}\langle \varphi,\varphi \rangle\right)\in\R
	\end{align*}
	is a characteristic function on $\mathcal{N}$. Using the Bocher-Minlos theorem (see, e.g., \cite{BK95} or \cite{Oba94}), the following definition from \cite{GJ16} makes sense.
	\begin{defn}
		For $0<\beta\leq 1$, the \textit{Mittag-Leffler measure} $\mu_\beta$ is defined as the unique probability measure on the space $(\mathcal{N}',\mathcal{C}_\sigma(\mathcal{N}'))$ whose characteristic function is
		\begin{equation*}
			\int_{\mathcal{N}'} e^{i\langle \omega,\varphi \rangle}\de\mu_\beta(\omega) = E_\beta\left(-\frac{1}{2}\langle \varphi,\varphi \rangle\right),\qquad \varphi\in \mathcal{N}.
		\end{equation*}
		The corresponding $L^p$ spaces of complex-valued functions are denoted by $L^p(\mu_\beta):=L^p(\mathcal{N}',\mu_\beta;\C)$ for $p\geq 1$ with corresponding norm $\|\cdot\|_{L^p(\mu_\beta)}$. For $p=2$, the corresponding scalar product is denoted by $(\!(\cdot,\cdot)\!)_{L^2(\mu_\beta)}$.
	\end{defn}
	
	\begin{rmk}
		The class of Mittag-Leffler measures on $\mathcal{N}'$ includes the following.
		\begin{itemize}
			\item For $\beta = 1$, the Mittag-Leffler measure $\mu_1$ on $\mathcal{N}'$ is the usual Gaussian measure on $\mathcal{N}'$ with covariance given by the scalar product in $\mathcal{H}$.
			\item If $\mathcal{H}=\mathcal{N}=\R^n$, $n\in\N$, the Mittag-Leffler measure on $\mathcal{N}'=\R^n$ is called the $n$-dimensional Mittag-Leffler measure, and is denoted by $\mu^n_\beta$. This has been studied in \cite{Sch92}.
			\item The measure $\mu_\beta$ on $\mathcal{S}'(\R)$ is also called the grey noise (reference) measure in \cite{GJRS15,GJ16}.
		\end{itemize}
	\end{rmk}
	In \cite{GJRS15,GJ16}, the following properties of the measure $\mu_\beta$ are obtained.
	\begin{prop}\label{p:bc01}
		For any $\varphi\in\mathcal{N}$ and $n\in\N_0$,
		\begin{gather}
			\int_{\mathcal{N}'} \langle \omega,\varphi\rangle^{2n+1}\de\mu_\beta(\omega) = 0;\notag\\
			\int_{\mathcal{N}'} \langle \omega,\varphi\rangle^{2n}\de\mu_\beta(\omega) = \frac{(2n)!}{2^n\Gamma(\beta n + 1)}\langle \varphi,\varphi \rangle^{n}.\notag
		\end{gather}
		In particular, for all $\varphi,\psi\in\mathcal{N}$,
		\begin{gather}\label{e:ml2}
			\|\langle\cdot,\varphi\rangle\|_{L^2(\mu_\beta)}^2 = \frac{1}{\Gamma(\beta+1)}|\varphi|_0^2,\\
			\int_{\mathcal{N}'} \langle \omega,\varphi\rangle \langle \omega,\psi\rangle\de\mu_\beta(\omega) = \frac{1}{\Gamma(\beta+1)}\langle \varphi,\psi\rangle. \notag
		\end{gather}
	\end{prop}
	
	Equation~(\ref{e:ml2}) allows us to define $\langle \cdot,\eta\rangle$, $\eta\in \mathcal{H}$, as an $L^2(\mu_\beta)$-limit of the sequence $(\langle \cdot,\varphi_n\rangle)_{n=1}^\infty$, where $(\varphi_n)_{n=1}^\infty$ is a sequence in $\mathcal{N}$ converging to $\eta$ in $\mathcal{H}$. Furthermore, this limit is independent of the sequence approximating $\eta$. It was shown in \cite{GJRS15,GJ16,Jah15} that Proposition~\ref{p:bc01} holds for elements in $\mathcal{H}$. Moreover, as a random variable on the probability space $(\mathcal{N}',\mathcal{C}_\sigma(\mathcal{N}'),\mu_\beta)$, its characteristic function is given by
	\begin{align}\label{e:chah}
		\int_{\mathcal{N}'}e^{ip\langle \cdot,\eta\rangle}\de\mu_\beta(\omega) = E_\beta\left(-\frac{p^2}{2}|\eta|_0^2\right),\quad p\in\R.
	\end{align}
	
	\begin{prop}\label{p:indmub}
		Let $\{\varphi_1,\dots,\varphi_n\}$, $n\in\mathbb{N}$, be an orthonormal set in $\mathcal{H}$. The random variables $\langle\cdot,\varphi_1\rangle$, $\dots$, $\langle\cdot,\varphi_n\rangle$ on the probability space $(\mathcal{N}',\mathcal{C}_\sigma(\mathcal{N}'),\mu_\beta)$ are independent if and only if $n=1$ or $\beta=1$.
	\end{prop}
	\pf Using (\ref{e:chah}), for all real numbers $p_1$, $\dots$, $p_n$,
	\begin{align*}
		\int_{\mathcal{N}'} \exp\left(i\sum_{r=1}^n p_r\langle\omega,\varphi_r\rangle\right)\,\text{\rm d}\mu_\beta(\omega) =  E_\beta\left(-\frac{1}{2}\sum_{r=1}^n p_r^2\right).
	\end{align*}
	If $n\neq 1$, then independence of $\langle\cdot,\varphi_1\rangle$, $\dots$, $\langle\cdot,\varphi_n\rangle$ holds if and and only if $E_\beta(-(x+y)) = E_\beta(-x)E_\beta(-y)$ for all $x,y\geq 0$. As $E_\beta$ is entire and $E_\beta(0)=1$, the identity theorem from complex analysis implies that this holds if and only if $E_\beta$ is the exponential function, that is, $\beta=1$. \qed
	
	\section{Mittag-Leffler analysis on product spaces}
	\subsection{Products of Mittag-Leffler spaces}
	Let $d\in\N$, and consider the real separable Hilbert space $L^2_d(\R):=L^2(\R,\!\de x;\R^d)$ of $\R^d$-valued square integrable functions on $\R$ with respect to the Lebesgue measure, which is isomorphic to the (external) direct sum $\bigoplus_{k=1}^d L^2(\R)$ of $d$ copies of $L^2(\R)$. In this case, every $\eta\in L^2_d(\R)$ can be written uniquely as
	\begin{align*}
		\eta = \sum_{k=1}^d \eta_k\mathbf{e}_k := (\eta_1,\dots,\eta_d),
	\end{align*}
	where $\eta_k\in L^2(\R)$ for each $k=1,\dots,d$, and $\{\mathbf{e}_1,\dots,\mathbf{e}_d\}$ is the canonical basis of $\R^d$. The norm of $\eta$ induced by the scalar product in $L^2_d(\R)$ is given by
	\begin{align*}
		|\eta|_{L^2_d(\R)}^2 = \sum_{k=1}^d |\eta_k|_{L^2(\R)}^2.
	\end{align*}
	 A dense subspace of $L^2_d(\R)$ is the space $\mathcal{S}_d(\R):=\bigoplus_{k=1}^d\mathcal{S}(\R)$, the (external) direct sum of $d$ copies of $\mathcal{S}(\R)$. This representation of $\mathcal{S}_d(\R)$ as a direct sum is useful to utilize some results in general duality theory; the details are as follows. If we equip $\mathcal{S}_d(\R)$ with the locally convex direct sum topology, then $\mathcal{S}_d(\R)$ is a nuclear space (see,~e.g.,~Proposition~50.1 in~\cite{Tre67}). In fact, if we let $\mathcal{H}_{d,p}:=\bigoplus_{k=1}^d \mathcal{H}_p$, where $(\mathcal{H}_p)_{p\in\N}$ are the Hilbert spaces defining the topology of $\mathcal{S}(\R)$, then $\mathcal{S}_d(\R)$ is the projective limit of the real separable Hilbert spaces $(\mathcal{H}_{d,p})_{p\in\N}$ satisfying (N1) and (N2). Here, the induced norm on $\mathcal{H}_{d,p}$ is given by
	\begin{align}\label{e:vsdnm}
		|\varphi|_p^2 := \sum_{k=1}^d |\varphi_k|_p^2,\quad \varphi\in\mathcal{H}_{d,p},\ p\in\N,
	\end{align}
	where the norm $|\cdot|_p$ appearing on the right-hand side of (\ref{e:vsdnm}) is the norm on $\mathcal{H}_p$. For notational convenience, we identify the norm $|\cdot|_0$ with the norm on $L^2_d(\R)$. Thus, $\mathcal{S}_d(\R)$ is a nuclear Frech\'{e}t space. Together with the dual space $\mathcal{S}_d'(\R)$, we obtain the nuclear triple
	\begin{align*}
		\mathcal{S}_d(\R) \subset L^2_d(\R) \subset \mathcal{S}_d'(\R).
	\end{align*}
	Recalling that inductive limit topology on $\mathcal{S}_d'(\R)$ coincides with the strong topology, $\mathcal{S}_d'(\R)$ is topologically isomorphic to the product space $\mathcal{S}'(\R)^d$ via the canonical identification
	\begin{align*}
		\mathcal{S}_d'(\R) \ni \omega \leftrightarrow (\omega_1,\dots,\omega_d) \in \mathcal{S}'(\R)^d,\quad \omega_k:= \omega\circ\iota_k,
	\end{align*}
	where for $k=1,\dots,d$,  $\iota_k$ is the canonical injection of the $k^\text{th}$ copy of $\mathcal{S}(\R)$ into $\mathcal{S}_d(\R)$ (see,~e.g., Proposition~14 (IV, p.~12) in~\cite{Bou87}, or 18.10~in~\cite{KN63}). Hence, if $\omega\in\mathcal{S}_d'(\R)$ and $\varphi\in\mathcal{S}_d(\R)$, then (see,~e.g.,~14.7~in~\cite{KN63})
	\begin{align*}
		\langle\omega,\varphi\rangle = \sum_{k=1}^d \langle\omega_k,\varphi_k\rangle,
	\end{align*}
	where the dual pairing on the right-hand side is that of between $\mathcal{S}'(\R)$ and $\mathcal{S}(\R)$. Moreover, recalling that  $\mathcal{C}_\sigma(\mathcal{S}_d'(\R))$ coincides with the Borel $\sigma$-algebra generated by the strong topology on $\mathcal{S}_d'(\R)$, it also coincides with the $d$-fold product $\sigma$-algebra $\mathcal{C}_\sigma(\mathcal{S}'(\R))^{\otimes d}$, since the map
	$$
		\mathcal{S}_d'(\R)\ni \omega \mapsto (\langle\omega,\varphi^1\rangle, \dots, \langle\omega,\varphi^n\rangle) = \sum_{k=1}^{d} (\langle\omega_k,\varphi^1_k\rangle, \dots, \langle\omega_k,\varphi^n_k\rangle) \in \R^n
	$$
	for $n\in\N$ and $\varphi^1,\dots,\varphi^n\in \mathcal{S}_d(\R)$ is $\mathcal{C}_\sigma(\mathcal{S}'(\R))^{\otimes d}$-$\mathcal{B}(\R^n)$ measurable.
	
	\bigskip
	
	Let $\varphi\in \mathcal{S}_d(\R)$. Throughout this paper, we consider the following $\mathcal{C}_\sigma(\mathcal{S}_d'(\R))$-$\mathcal{B}(\R^n)$ measurable map $G(\cdot,\varphi)$ on $\mathcal{S}_d'(\R)$ defined by
	\begin{align*}
		G(\omega,\varphi):=(\langle\omega_1,\varphi_1\rangle,\dots,\langle\omega_d,\varphi_d\rangle)\in\R^d,\quad \omega\in\mathcal{S}_d'(\R),
	\end{align*}
	and we would like to have a probability measure on $(\mathcal{S}_d'(\R),\mathcal{C}_\sigma(\mathcal{S}_d'(\R)))$ such that $G(\cdot,\varphi)$ is a random vector whose components are mutually independent random variables with respect to Mittag-Leffler measure $\mu_\beta$ on $\mathcal{S}'(\R)$. A natural way to obtain this property is to use the $d$-fold product measure of the Mittag-Leffler measure, denoted by $\mu_{\beta}^{\otimes d}$, whose characteristic function is given by
	\begin{align*}
		\int_{\mathcal{S}_d'(\R)} e^{i\langle\omega,\varphi\rangle}\de\mu_{\beta}^{\otimes d}(\omega) = \prod_{k=1}^{d}E_\beta\left(-\frac{1}{2}\langle\varphi_k,\varphi_k\rangle\right),\qquad \varphi\in \mathcal{S}_d(\R).
	\end{align*}
	For $p\geq 1$, the $L^p$-space of complex-valued functions on $\mathcal{S}_d'(\R)$ is denoted by $L^p(\mu_{\beta}^{\otimes d}):=L^p(\mathcal{S}_d'(\R),\mu_{\beta}^{\otimes d};\C)$ with corresponding norm $\|\cdot\|_{L^p(\mu_{\beta}^{\otimes d})}$. For $p=2$, the corresponding scalar product, denoted by $(\!(\cdot,\cdot)\!)_{L^2(\mu_{\beta}^{\otimes d})}$, is defined as follows:
	\begin{align*}
		(\!(F,G)\!)_{L^2(\mu_{\beta}^{\otimes d})} := \int_{\mathcal{S}_d'(\R)} F(\omega)\overline{G(\omega)}\de\mu_{\beta}^{\otimes d}(\omega),\quad F,G\in L^2(\mu_{\beta}^{\otimes d}).
	\end{align*}
	
	\begin{prop}\label{p:mmnl}
		For any $\varphi\in\mathcal{S}_d(\R)$ and $n\in\N_0$,
		\begin{gather}
			\int_{\mathcal{S}_d'(\R)} \langle\omega,\varphi\rangle^{2n+1}\de\mu_{\beta}^{\otimes d}(\omega) = 0;\label{e:mmnl1}\\
			\int_{\mathcal{S}_d'(\R)} \langle\omega,\varphi\rangle^{2n}\de\mu_{\beta}^{\otimes d}(\omega) = \dfrac{(2n)!}{2^n} \sum_{r} \dfrac{\langle\varphi_1,\varphi_1\rangle^{r_1}\cdots\langle\varphi_d,\varphi_d\rangle^{r_d}}{\Gamma(\beta r_1 + 1)\cdots\Gamma(\beta r_d + 1)}, \label{e:mmnl2}
		\end{gather}
		where the sum in (\ref{e:mmnl2}) is taken over all $r:=(r_1,\dots,r_d)\in\N_0^d$ such that $r_1+\cdots+r_d = n$, and we use the convention that $\langle\varphi_k,\varphi_k\rangle^0=1$, even if $\varphi_k=0$. In particular, for all $\varphi,\psi\in\mathcal{S}_d(\R)$,
		\begin{gather}\label{e:l2nl}
			\|\langle\cdot,\varphi\rangle\|_{L^2(\mu_{\beta}^{\otimes d})}^2 = \frac{1}{\Gamma(\beta+1)}|\varphi|_0^2,\\
			\int_{\mathcal{S}_d'(\R)} \langle\omega,\varphi\rangle\langle\omega,\psi\rangle\de\mu_{\beta}^{\otimes d}(\omega) = \frac{1}{\Gamma(\beta+1)}\langle\varphi,\psi\rangle.\notag
		\end{gather}
	\end{prop}
	
	\pf The multinomial theorem yields that for $m\in\N_0$,
	\begin{align*}
		\int_{\mathcal{S}_d'(\R)} \langle\omega,\varphi\rangle^{m}\de\mu_{\beta}^{\otimes d}(\omega) &= \int_{\mathcal{S}_d'(\R)} \left(\sum_{k=1}^d \langle\omega_k,\varphi_k\rangle\right)^{m}\text{\rm d}\mu_{\beta}^{\otimes d}(\omega) = \sum_{p} \dfrac{m!}{p_1!\cdots p_d!}\prod_{k=1}^d \int_{\mathcal{S}'(\R)} \langle\omega_k,\varphi_k\rangle^{p_k}\de\mu_\beta(\omega_k),
	\end{align*}
	where the sum is taken over all $p:=(p_1,\dots,p_d)\in\N_0^d$ such that $p_1+\cdots+p_d = m$. Equations (\ref{e:mmnl1})-(\ref{e:mmnl2}) then follow directly from Proposition~\ref{p:bc01}. \qed
	
	\begin{rmk}\label{r:l2rv}
		Equation~(\ref{e:l2nl}) allows us to define $\langle\cdot,\eta\rangle$ for $\eta\in L^2_d(\R)$ as an $L^2(\mu_{\beta}^{\otimes d})$-limit of the sequence $(\langle\cdot,\varphi_n\rangle)_{n\in\N}$, where $(\varphi_n)_{n\in\N}$ is a sequence in $\mathcal{S}_d(\R)$ that converges to $\eta$ with respect to the $L^2_d(\R)$ norm. Using the proofs similar to that of the case for the one-dimensional Mittag-Leffler measure, this limit is independent of the approximating sequence $(\varphi_n)_{n\in\N}$ of $\eta$, and that Proposition~\ref{p:mmnl} holds for elements in $L^2_d(\R)$. Moreover, for $\mu_{\beta}^{\otimes d}$-almost all $\omega\in\mathcal{S}_d'(\R)$,
		\begin{align*}
			\langle\omega,\eta\rangle = \sum_{k=1}^d \langle\omega,\eta_k\mathbf{e}_k\rangle = \sum_{k=1}^{d}\langle\omega_k,\eta_k\rangle.
		\end{align*}
	\end{rmk}
	
	\subsection{The Appell system}
	The measures $\mu_{\beta}^{\otimes d}$ are generally non-Gaussian, and so in the same manner as that of $\mu_\beta$, the construction of test and distribution spaces on $\mu_{\beta}^{\otimes d}$ uses the Appell system, introduced by \cite{KSWY98}.
	
	\subsubsection{Compatibility}
	We refer to \cite{GJRS15} for the requirements to use the Appell system: given a nuclear triple $\mathcal{N}\subset\mathcal{H}\subset\mathcal{N}'$, any measure $\mu$ on $(\mathcal{N}',\mathcal{C}_\sigma(\mathcal{N}'))$ must satisfy the following properties:
	\begin{itemize}
		\item[(A1)] The measure $\mu$ has an analytic Laplace transform in a neighborhood of zero, that is, the map
		\begin{align*}
			\mathcal{N}_\C\ni \varphi\mapsto l_\mu(\varphi):= \int_{\mathcal{N}'} e^{\langle\omega,\varphi\rangle}\de\mu(\omega)\in\C
		\end{align*}
		is holomorphic on a neighborhood $\mathcal{U}_0\subset\mathcal{N}_\C$ of zero.
		\item[(A2)] For any nonempty open subset $\mathcal{U}\subset\mathcal{N}'$, we have $\mu(\mathcal{U})>0$.
	\end{itemize}
	
	In the following, we show that $\mu_{\beta}^{\otimes d}$ satisfies (A1) and (A2).
	
	\begin{lem}\label{p:exp}
		Let $\varphi\in\mathcal{S}_d(\R)$ and $\lambda\in\R$. Then the exponential map $\mathcal{S}_d'(\R)\ni\omega\mapsto e^{|\lambda\langle\omega,\varphi\rangle|}$ is integrable with respect to $\mu_{\beta}^{\otimes d}$ and
		\begin{align*}
			\int_{\mathcal{S}_d'(\R)} e^{\lambda\langle\omega,\varphi\rangle}\de\mu_{\beta}^{\otimes d}(\omega) = \prod_{k=1}^d E_\beta\left(\dfrac{\lambda^2}{2}\langle\varphi_k,\varphi_k\rangle\right).
		\end{align*}
	\end{lem}
	\pf By Lemma~4.1 in \cite{GJRS15},
	\begin{align*}
		\int_{\mathcal{S}_d'(\R)} e^{|\lambda\langle\omega,\varphi\rangle|}\de\mu_{\beta}^{\otimes d}(\omega) &\leq \int_{\mathcal{S}_d'(\R)} \exp\left(\sum_{k=1}^d \Big|\lambda\langle\omega_k,\varphi_k\rangle\Big|\right)\de\mu_{\beta}^{\otimes d}(\omega) = \prod_{k=1}^d \int_{\mathcal{S}'(\R)} e^{\big|\lambda \langle\omega_k,\varphi_k\rangle\big|}\de\mu_\beta(\omega_k) < \infty,
	\end{align*}
	and
	\begin{align*}
		\int_{\mathcal{S}_d'(\R)} e^{\lambda\langle\omega,\varphi\rangle}\de\mu_{\beta}^{\otimes d}(\omega) = \prod_{k=1}^d \int_{\mathcal{S}'(\R)} e^{\lambda \langle\omega_k,\varphi_k\rangle}\de\mu_\beta(\omega_k) = \prod_{k=1}^d E_\beta\left(\dfrac{\lambda^2}{2}\langle\varphi_k,\varphi_k\rangle\right). \tag*{\qed}
	\end{align*}
	
	\begin{prop}\label{p:a1}
		The map
		\begin{align*}
			\mathcal{S}_d(\R)_\C\ni \varphi \mapsto l_{\mu_{\beta}^{\otimes d}}(\varphi):= \int_{\mathcal{S}_d'(\R)} e^{\langle\omega,\varphi\rangle}\de\mu_{\beta}^{\otimes d}(\omega)\in\C
		\end{align*}
		is a holomorphic map from $\mathcal{S}_d(\R)_\C$ to $\C$.
	\end{prop}
	\pf Note that $l_{\mu_{\beta}^{\otimes d}}$ is locally bounded on $\mathcal{S}_d(\R)_\C$. Indeed, if $\varphi:=\varphi^1 + i\varphi^2\in\mathcal{S}_d(\R)_\C$, then Lemma~\ref{p:exp} implies that
	\begin{align*}
		|l_{\mu_{\beta}^{\otimes d}}(\varphi)| \leq \int_{\mathcal{S}_d'(\R)} e^{\langle\omega,\varphi^1\rangle}\de\mu_{\beta}^{\otimes d}(\omega) = \prod_{k=1}^d E_\beta\left(\dfrac{1}{2}\langle\varphi^1_k,\varphi^1_k\rangle\right) <\infty.
	\end{align*}
	Now we show that $l_{\mu_{\beta}^{\otimes d}}$ is G-holomorphic, that is, the map $\C\ni z\mapsto f(z):=l_{\mu_{\beta}^{\otimes d}}(\varphi^0 + z\varphi)$, where $\varphi^0,\varphi\in \mathcal{S}_d(\R)_\C$, is holomorphic on some neighborhood of zero in $\C$. Note that $f$ is continuous: given $z\in\C$ and a sequence $(z_n)_{n\in\N}$ in $\C$ converging to $z$, the following estimate holds for sufficiently large $n$:
	\begin{align*}
		|\exp(\langle\omega,\varphi^0 + z_n\varphi\rangle)| \leq \exp(|\langle\omega,\varphi^0\rangle|) \exp\big((1+|z|)|\langle\omega,\varphi\rangle|\big),
	\end{align*}
	and thus continuity of $f$ follows from Lemma~\ref{p:exp}, Cauchy-Schwarz inequality, and Lebesgue dominated convergence theorem. Moreover, if $\gamma$ is a closed, bounded curve in $\C$, then the compactness of $\gamma$ allows us to use Fubini's theorem:
	\begin{align*}
		\int_{\gamma}\int_{\mathcal{S}_d'(\R)} \exp(\langle\omega,\varphi^0 + z\varphi\rangle)\de\mu_{\beta}^{\otimes d}(\omega)\de z = \int_{\mathcal{S}_d'(\R)}\int_{\gamma} \exp(\langle\omega,\varphi^0 + z\varphi\rangle)\de z\de\mu_{\beta}^{\otimes d}(\omega) = 0,
	\end{align*}
	where the last equality holds as the exponential function is holomorphic on $\C$. By Morera's theorem, $f$ is holomorphic on $\C$, and thus, $l_{\mu_{\beta}^{\otimes d}}$ is G-holomorphic. This implies that $l_{\mu_{\beta}^{\otimes d}}$ is holomorphic (see~\cite{Din81}). \qed
	
	\bigskip
	
	Using Lemma~\ref{p:exp}, the identity theorem from complex analysis to the function $f$ in the proof of Proposition~\ref{p:a1}, and the density of $\mathcal{S}_d(\R)_\C$ in $L^2_d(\R)_\C$, the following corollary holds.
	
	\begin{cor}\label{c:excom}
		For $\eta\in L^2_d(\R)_\C$ and $z\in\C$, the map $\mathcal{S}_d'(\R)\ni\omega\mapsto e^{|z\langle\omega,\eta\rangle|}$ is integrable with respect to $\mu_{\beta}^{\otimes d}$ and
		\begin{align*}
			\int_{\mathcal{S}_d'(\R)} e^{z\langle\omega,\eta\rangle}\de\mu_{\beta}^{\otimes d}(\omega) = \prod_{k=1}^d E_\beta\left(\dfrac{z^2}{2}\langle\eta_k,\eta_k\rangle\right).
		\end{align*}
	\end{cor}
	
	\bigskip
	
	The proof that $\mu_{\beta}^{\otimes d}$ satisfies (A2) is straightforward.
	
	\begin{prop}\label{p:a2}
		For any nonempty open subset $\mathcal{U}\subset\mathcal{S}_d'(\R)$, we have $\mu_{\beta}^{\otimes d}(\mathcal{U})>0$.
	\end{prop}
	\pf Let $\mathcal{U}$ be a nonempty open subset of $\mathcal{S}_d'(\R)$. As $\mathcal{S}_d'(\R)$ is a finite product space, there exist nonempty open sets $\mathcal{U}_1,\dots,\mathcal{U}_d$ of $\mathcal{S}'(\R)$ such that $\mathcal{U}_1 \times \cdots \times \mathcal{U}_d \subset \mathcal{U}$. Since the measure $\mu_{\beta}$ on $\mathcal{S}'(\R)$ satisfies (A2) by Theorem~4.5 in~\cite{GJRS15}, we have
	$$
		\mu_{\beta}^{\otimes d}(\mathcal{U}) \geq \mu_{\beta}(\mathcal{U}_1)\cdot\ldots\cdot\mu_{\beta}(\mathcal{U}_d) > 0. \eqno{\qed}
	$$
	
	\subsubsection{Construction of test functions and distributions}
	
	Now, we proceed to construct the test function space and the distribution space via the Appell system. Details of the construction and most of the notations and statements of this subsection can be found in \cite{GJRS15} and references therein.
	
	First, we introduce the $\mu_{\beta}^{\otimes d}$-exponential by
	\begin{align*}
		e_{\mu_{\beta}^{\otimes d}}(\varphi;\omega):=\dfrac{e^{\langle \omega,\varphi\rangle}}{l_{\mu_{\beta}^{\otimes d}}(\varphi)},\quad \varphi\in\mathcal{S}_d(\R)_\C,\ \omega\in\mathcal{S}'_d(\R)_\C.
	\end{align*}
	Since $l_{\mu_{\beta}^{\otimes d}}(0)=1$ and $l_{\mu_{\beta}^{\otimes d}}$ is holomorphic, there exists a neighborhood $\mathcal{U}_0\subset\mathcal{S}_d(\R)_\C$ of zero such that $l_{\mu_{\beta}^{\otimes d}}(\varphi)\neq 0$ for all $\varphi\in\mathcal{U}_0$, and thus, the $\mu_{\beta}^{\otimes d}$-exponential is well-defined on $\mathcal{U}_0$. In this case, the $\mu_{\beta}^{\otimes d}$-exponential can be expressed as a power series, i.e.,
	\begin{align*}
		e_{\mu_{\beta}^{\otimes d}}(\varphi;\omega) = \sum_{n=0}^{\infty} \dfrac{1}{n!}\langle P_n^{\mu_{\beta,d}}(\omega), \varphi^{\otimes n}\rangle,\quad \varphi\in\mathcal{U}_0,\ \omega\in\mathcal{S}'_d(\R)_\C,
	\end{align*}
	for suitable mappings $P_n^{\mu_{\beta,d}}: \mathcal{S}'_d(\R)_\C \to (\mathcal{S}_d(\R)_\C^{\widehat{\otimes} n})'$. Using these mappings, every $\varphi\in\mathcal{P}(\mathcal{S}_d'(\R))$, the space of smooth polynomials on $\mathcal{S}_d'(\R)$, has the following unique representation:
	\begin{align}\label{e:pmu}
		\varphi = \sum_{n=0}^N\langle P_n^{\mu_{\beta,d}}(\cdot),\varphi^{(n)}\rangle,
	\end{align}
	for suitable $N\in\N_0$ and $\varphi^{(n)}\in\mathcal{S}_d(\R)_\C^{\widehat{\otimes} n}$. Moreover, every $\Phi\in\mathcal{P}_{\mu_{\beta}^{\otimes d}}'(\mathcal{S}_d'(\R))$, the dual space of $\mathcal{P}(\mathcal{S}_d'(\R))$ with respect to $L^2(\mu_{\beta}^{\otimes d})$, also has a unique representation:
	\begin{align}\label{e:qmu}
		\Phi=\sum_{n=0}^\infty Q_n^{\mu_{\beta,d}}(\Phi^{(n)}),
	\end{align}
	for suitable $\Phi^{(n)}\in (\mathcal{S}_d(\R)_\C^{\widehat{\otimes} n})'$. In this representation, $Q_n^{\mu_{\beta,d}}(\Phi^{(n)}):=D(\Phi^{(n)})^*\mathbf{1}$, where $D(\Phi^{(n)})^*$ is the adjoint of the continuous linear operator $D(\Phi^{(n)})$ on $\mathcal{P}(\mathcal{S}'(\R))$ defined on the monomials $\langle \cdot^{\otimes m},\varphi^{(m)}\rangle$ by
	\begin{align*}
		D(\Phi^{(n)})\langle \omega^{\otimes m},\varphi^{(m)}\rangle =
		\begin{cases}
			\dfrac{m!}{(m-n)!} \langle \omega^{\otimes (m-n)} \widehat{\otimes} \Phi^{(n)},\varphi^{(m)}\rangle, & m\geq n,\\
			0, & m<n,
		\end{cases}
	\end{align*}
	while $\mathbf{1}\in L^2(\mu_{\beta}^{\otimes d})$ is defined by $\mathbf{1}(\omega)\equiv 1$ for all $\omega\in\mathcal{S}_d'(\R)$. The Appell system $\mathbb{A}^{\mu_{\beta,d}}$ is then defined by the pair $(\mathbb{P}^{\mu_{\beta,d}},\mathbb{Q}^{\mu_{\beta,d}})$, where
	$$
		\mathbb{P}^{\mu_{\beta,d}} := \{\langle P_n^{\mu_{\beta,d}}(\cdot),\varphi^{(n)}\rangle: \varphi^{(n)}\in\mathcal{S}_d(\R)_\C^{\widehat{\otimes} n}, n\in\N_0\},\quad \mathbb{Q}^{\mu_{\beta,d}} := \{Q_n^{\mu_{\beta,d}}(\Phi^{(n)}) : \Phi^{(n)}\in(\mathcal{S}_d(\R)_\C^{\widehat{\otimes} n})', n\in\N_0\}.
	$$
	
	\bigskip
	
	The central property of the Appell system is the following bi-orthogonality relation:
	\begin{thm}{\rm\bf(Theorem~4.17 in \cite{KSWY98})}
		For all $\Phi^{(n)}\in(\mathcal{S}_d(\R)_\C^{\widehat{\otimes} n})'$ and $\varphi^{(m)}\in\mathcal{S}_d(\R)_\C^{\widehat{\otimes} m}$
		\begin{align}\label{t:biorth}
			\llangle Q_n^{\mu_{\beta,d}}(\Phi^{(n)}),\langle P_m^{\mu_{\beta,d}}(\cdot),\varphi^{(m)}\rangle\rrangle_{\mu_{\beta}} = \delta_{m,n}n!\langle \Phi^{(n)}, \varphi^{(n)}\rangle.
		\end{align}
	\end{thm}
	\noindent From here on, we always use the representations (\ref{e:pmu}) and (\ref{e:qmu}) for $\varphi\in\mathcal{P}(\mathcal{S}_d'(\R))$ and $\Phi\in\mathcal{P}_{\mu_{\beta}^{\otimes d}}'(\mathcal{S}_d'(\R))$, respectively.
	
	\bigskip
	
	For $p,q\in\N_0$, the space $(\mathcal{H}_{d,p})^1_{q,\mu_{\beta}^{\otimes d}}\subset L^2(\mu_{\beta}^{\otimes d})$ is defined as the completion of the space $\mathcal{P}(\mathcal{S}_d'(\R))$ with respect to the norm
	\begin{align*}
		\|\varphi\|_{p,q,\mu_{\beta}^{\otimes d}}^2 := \sum_{n=0}^N (n!)^2 2^{nq}|\varphi^{(n)}|_p^2,\quad\text{for }\varphi\in \mathcal{P}(\mathcal{S}_d'(\R)),
	\end{align*}
	while its dual $(\mathcal{H}_{d,-p})^{-1}_{-q,\mu_{\beta}^{\otimes d}}$ is the set of all $\Phi\in\mathcal{P}_{\mu_{\beta}^{\otimes d}}'(\mathcal{S}_d'(\R))$ such that the norm
	\begin{align*}
		\|\Phi\|_{-p,-q,\mu_{\beta}^{\otimes d}}^2 := \sum_{n=0}^\infty 2^{-qn}|\Phi^{(n)}|_{-p}^2 < \infty.
	\end{align*}
	The test function space $(\mathcal{S}_d(\R))^1_{\mu_{\beta}^{\otimes d}}$ is then defined as
	\begin{align*}
		(\mathcal{S}_d(\R))^1_{\mu_{\beta}^{\otimes d}} := \prlim_{p,q\in\N} (\mathcal{H}_{d,p})^1_{q,\mu_{\beta}^{\otimes d}},
	\end{align*}
	while the distribution space $(\mathcal{S}_d(\R))^{-1}_{\mu_{\beta}^{\otimes d}}$ is defined as
	\begin{align*}
		(\mathcal{S}_d(\R))^{-1}_{\mu_{\beta}^{\otimes d}} := \indlim_{p,q\in\N} (\mathcal{H}_{d,-p})^{-1}_{-q,\mu_{\beta}^{\otimes d}}.
	\end{align*}
	Hence, we obtain the following chain of continuous dense embeddings:
	\begin{align*}
		(\mathcal{S}_d(\R))^1_{\mu_{\beta}^{\otimes d}} \subset (\mathcal{H}_{d,p})^1_{q,\mu_{\beta}^{\otimes d}} \subset L^2(\mu_{\beta}^{\otimes d}) \subset (\mathcal{H}_{d,-p})^{-1}_{-q,\mu_{\beta}^{\otimes d}} \subset (\mathcal{S}_d(\R))^{-1}_{\mu_{\beta}^{\otimes d}}.
	\end{align*}
	The dual pairing between $(\mathcal{S}_d(\R))^{-1}_{\mu_{\beta}^{\otimes d}}$ and $(\mathcal{S}_d(\R))^1_{\mu_{\beta}^{\otimes d}}$ is denoted by $\llangle\cdot,\cdot\rrangle_{\mu_{\beta}^{\otimes d}}$, and is a bilinear extension of the scalar product on $L^2(\mu_{\beta}^{\otimes d})$ given by
	\begin{align*}
		\llangle F,\varphi\rrangle_{\mu_{\beta}^{\otimes d}} = (\!( F,\overline{\varphi} )\!)_{L^2(\mu_{\beta}^{\otimes d})},\quad \varphi\in\mathcal{P}(\mathcal{S}_d'(\R)),\ F\in L^2(\mu_{\beta}^{\otimes d}).
	\end{align*}
	The bi-orthogonality relation (\ref{t:biorth}) implies that for any $\varphi\in(\mathcal{S}_d(\R))^1_{\mu_{\beta}^{\otimes d}}$ and $\Phi\in(\mathcal{S}_d(\R))^{-1}_{\mu_{\beta}^{\otimes d}}$, we have
	\begin{align*}
		\llangle\Phi,\varphi\rrangle_{\mu_{\beta}^{\otimes d}} = \sum_{n=0}^\infty n! \langle \Phi^{(n)},\varphi^{(n)}\rangle.
	\end{align*}
	The same notation is used for the dual pairing between $(\mathcal{H}_{d,-p})^{-1}_{-q,\mu_{\beta}^{\otimes d}}$ and $(\mathcal{H}_{d,p})^1_{q,\mu_{\beta}^{\otimes d}}$. The set
	\begin{align*}
		\left\{ e_{\mu_{\beta}^{\otimes d}}(\varphi;\cdot):\ \varphi\in \mathcal{S}_d(\R)_\C,\ 2^q|\varphi|_p<1\right\}
	\end{align*}
	is total in $(\mathcal{H}_{d,p})^1_{q,\mu_{\beta}^{\otimes d}}$ and for any $\varphi$ belonging to the set $\mathcal{U}_{p,q}:=\{\varphi\in\mathcal{S}_d(\R)_\C:2^q|\varphi|_p<1\}$, we have $\|e_{\mu_{\beta}^{\otimes d}}(\varphi;\cdot)\|_{p,q,\mu_{\beta}^{\otimes d}}<\infty$.
	
	\subsubsection{Integral transforms and characterization theorems}
	
	As in the case of Gaussian analysis, the $S_{\mu_{\beta}^{\otimes d}}$-transform and the $T_{\mu_{\beta}^{\otimes d}}$-transform are defined as follows. For any $\Phi\in(\mathcal{S}_d(\R))^{-1}_{\mu_{\beta}^{\otimes d}}$ and $\varphi\in\mathcal{U}_0\subset \mathcal{S}_d(\R)_\C$, where $\mathcal{U}_0$ is a suitable neighborhood of zero, we define
	\begin{align*}
		S_{\mu_{\beta}^{\otimes d}}\Phi(\varphi) := \llangle\Phi,e_{\mu_{\beta}^{\otimes d}}(\varphi;\cdot)\rrangle_{\mu_{\beta}^{\otimes d}},\\
		T_{\mu_{\beta}^{\otimes d}}\Phi(\varphi) := \llangle\Phi,e^{i\langle\cdot,\varphi\rangle}\rrangle_{\mu_{\beta}^{\otimes d}}.
	\end{align*}
	For a vector $\boldsymbol{\Phi}$ with components in $(\mathcal{S}_d(\R))^{-1}_{\mu_{\beta}^{\otimes d}}$, its $S_{\mu_{\beta}^{\otimes d}}$-transform and $T_{\mu_{\beta}^{\otimes d}}$-transform is a vector field defined on $\mathcal{U}_0$ whose components are the $S_{\mu_{\beta}^{\otimes d}}$-transform and $T_{\mu_{\beta}^{\otimes d}}$-transform, respectively, of the corresponding components of $\boldsymbol{\Phi}$. Of course, we have the following relationships between the two transforms:
	\begin{align}\label{e:st}
		T_{\mu_{\beta}^{\otimes d}}\Phi(\varphi) = l_{\mu_{\beta}^{\otimes d}}(i\varphi)\cdot S_{\mu_{\beta}^{\otimes d}}\Phi(i\varphi) = \prod_{k=1}^d E_\beta\left(-\dfrac{1}{2}\langle\varphi_k,\varphi_k\rangle\right)\cdot S_{\mu_{\beta}^{\otimes d}}\Phi(i\varphi),\quad\text{for } \varphi\in\mathcal{U}_0.
	\end{align}
	The characterization theorem for the space $(\mathcal{S}_d(\R))^{-1}_{\mu_{\beta}^{\otimes d}}$ via the $S_{\mu_{\beta}^{\otimes d}}$-transform is done using the spaces of holomorphic functions on $\mathcal{S}_d(\R)_\C$. We denote by $\text{Hol}_0\left(\mathcal{S}_d(\R)_\C\right)$ the space of holomorphic functions at zero, where we identify two functions which coincides on a neighborhood of zero. See \cite{KSWY98} for the details and proof of the following characterization theorem.
	\begin{thm}{\rm\bf(Theorem~8.34 in \cite{KSWY98})}\label{t:chrkd}
		The $S_{\mu_{\beta}^{\otimes d}}$-transform is a topological isomorphism from $(\mathcal{S}_d(\R))^{-1}_{\mu_{\beta}^{\otimes d}}$ to $\text{\rm Hol}_0\left(\mathcal{S}_d(\R)_\C\right)$.
	\end{thm}
	As in the case of $\mu_\beta$, a corollary of the characterization theorem is the result characterizing integrable maps with values in $(\mathcal{S}_d(\R))^{-1}_{\mu_{\beta}^{\otimes d}}$ in a weak sense. The proof is similar to that of the characterization theorem for $(\mathcal{S}(\R))^{-1}_{\mu_{\beta}}$; see Theorem~4.10 in \cite{GJRS15} or Theorem~2.2.2 in \cite{Jah15}.
	\begin{thm}\label{t:chrint}
		Let $(T,\mathcal{B},\nu)$ be a measure space and $\Phi_t\in(\mathcal{S}_d(\R))^{-1}_{\mu_{\beta}^{\otimes d}}$ for all $t\in T$. Let $\mathcal{U}_0\subset\mathcal{S}_d(\R)_\C$ be a suitable neighborhood of zero and $C<\infty$ such that
		\begin{itemize}
			\item[(i)] the map $T\ni t\mapsto S_{\mu_{\beta}^{\otimes d}}\Phi_t(\varphi)\in\C$ is measurable for all $\varphi\in\mathcal{U}_0$; and
			\item[(ii)] $\int_{T}|S_{\mu_{\beta}^{\otimes d}}\Phi_t(\varphi)|\de\nu(t)\leq C$ for all $\varphi\in\mathcal{U}_0$.
		\end{itemize}
		Then there exists a unique $\Psi\in(\mathcal{S}_d(\R))^{-1}_{\mu_{\beta}^{\otimes d}}$ such that for all $\varphi\in\mathcal{U}_0$,
		\begin{align*}
			S_{\mu_{\beta}^{\otimes d}}\Psi(\varphi) = \int_{T} S_{\mu_{\beta}^{\otimes d}}\Phi_t(\varphi)\de\nu(t).
		\end{align*}
		We denote $\Psi$ by $\int_{T} \Phi_t\de\nu(t)$ and call it the weak integral of $\Phi:=(\Phi_t)_{t\in T}$.
	\end{thm}
	
	Another corollary of the characterization theorem is the result characterizing strongly convergent sequences in $(\mathcal{S}_d(\R))^{-1}_{\mu_{\beta}^{\otimes d}}$ in terms of their $S_{\mu_{\beta}^{\otimes d}}$-transform. Again, the proof follows similarly as in the corresponding result for $(\mathcal{S}(\R))^{-1}_{\mu_{\beta}}$; see Theorem~2.12 in \cite{GJ16} or Theorem~2.3.1 in \cite{Jah15}.
	\begin{thm}\label{t:chrseq}
		A sequence $(\Phi_n)_{n\in\N}$ in $(\mathcal{S}_d(\R))^{-1}_{\mu_{\beta}^{\otimes d}}$ converges strongly in $(\mathcal{S}_d(\R))^{-1}_{\mu_{\beta}^{\otimes d}}$ if and only if there exist $p,q\in\N$ such that
		\begin{itemize}
			\item[(i)] $(S_{\mu_{\beta}^{\otimes d}}\Phi_n(\varphi))_{n\in\N}$ is a Cauchy sequence for all $\varphi\in\mathcal{U}_{p,q}$;
			\item[(ii)] for each $n\in\N$, $S_{\mu_{\beta}^{\otimes d}}\Phi_n$ is holomorphic on $\mathcal{U}_{p,q}$, and there exists a constant $C<\infty$ such that $|S_{\mu_{\beta}^{\otimes d}}\Phi_n(\varphi)|\leq C$ for all $\varphi\in\mathcal{U}_{p,q}$ and $n\in\N$.
		\end{itemize}
	\end{thm}
	
	A consequence of Theorem~\ref{t:chrseq} that is used for applications to stochastic differential equations is the following sufficient condition for the derivative and the $S_{\mu_{\beta}^{\otimes d}}$-transform to commute.
	
	\begin{cor}\label{c:dvst}
		Let $I\subset\R$ be an interval and $\Phi:I\to(\mathcal{S}_d(\R))^{-1}_{\mu_{\beta}^{\otimes d}}$. Assume that there exist $p,q\in\N$ such that
		\begin{enumerate}
			\item[(i)] for all $s\in I$, $S_{\mu_{\beta}^{\otimes d}}\Phi(s)$ is holomorphic on $\mathcal{U}_{p,q}$;
			\item[(ii)] for each $\varphi\in\mathcal{U}_{p,q}$, the map $I\ni t\mapsto S_{\mu_{\beta}^{\otimes d}}\Phi(t)(\varphi)\in\C$
			is differentiable;
			\item[(iii)] there exists a constant $C<\infty$ such that
			\begin{align*}
				\left|\frac{\text{\rm d}}{\text{\rm d}t}S_{\mu_{\beta}^{\otimes d}}\Phi(t)(\varphi)\right| \leq C,\quad\text{ for all } t\in I,\ \varphi\in\mathcal{U}_{p,q}.
			\end{align*}
		\end{enumerate}
		Then $\Phi$ is differentiable in $(\mathcal{S}_d(\R))^{-1}_{\mu_{\beta}^{\otimes d}}$ at all $t\in I$, that is, for each $t\in I$,
		\begin{align*}
			\frac{\text{\rm d}}{\text{\rm d}t}\Phi(t) := \lim_{h\to 0}\frac{1}{h}\left(\Phi(t+h)-\Phi(t)\right)
		\end{align*}
		exists as an element of $(\mathcal{S}_d(\R))^{-1}_{\mu_{\beta}^{\otimes d}}$. Moreover, for each $\varphi\in\mathcal{U}_{p,q}$ and $t\in I$, 
		\begin{align*}
			S_{\mu_{\beta}^{\otimes d}}\frac{\text{\rm d}}{\text{\rm d}t}\Phi(t)(\varphi) = \frac{\text{\rm d}}{\text{\rm d}t}S_{\mu_{\beta}^{\otimes d}}\Phi(t)(\varphi).
		\end{align*}
	\end{cor}
	\pf Let $t\in I$ and $(h_n)_{n\in\N}$ be a nonzero sequence in $\R$ with $t+h_n\in I$ and $h_n\to 0$ as $n\to\infty$. For each $n\in\N$, set
	\begin{align*}
		\Psi_n := \frac{1}{h_n}(\Phi(t+h_n)-\Phi(t))\in (\mathcal{S}_d(\R))^{-1}_{\mu_{\beta}^{\otimes d}}.
	\end{align*}
	Then for all $\varphi\in\mathcal{U}_{p,q}$, $(S_{\mu_{\beta}^{\otimes d}}\Psi_n(\varphi))_{n\in\N}$ is a Cauchy sequence, since
	\begin{align*}
		S_{\mu_{\beta}^{\otimes d}}\Psi_n(\varphi) = \frac{1}{h_n}(S_{\mu_{\beta}^{\otimes d}}\Phi(t+h_n)(\varphi)-S_{\mu_{\beta}^{\otimes d}}\Phi(t)(\varphi))\to \frac{\text{\rm d}}{\text{\rm d}t}S_{\mu_{\beta}^{\otimes d}}\Phi(t)(\varphi)\quad\text{as } n\to\infty.
	\end{align*}
	Moreover, we infer by the fundamental theorem of calculus that
	\begin{align*}
		|S_{\mu_{\beta}^{\otimes d}}\Psi_n(\varphi)| \leq \frac{1}{|h_n|}\left|\int_{t}^{t+h_n}\frac{\text{\rm d}}{\text{\rm d}s}S_{\mu_{\beta}^{\otimes d}}\Phi(s)(\varphi)\de s\right|\leq C.
	\end{align*}
	Thus, the sequence $(\Psi_n)_{n\in\N}$ fulfills the assumptions of Theorem~\ref{t:chrseq}, implying the existence of $\frac{\text{\rm d}}{\text{\rm d}t}\Phi(t)$ in $(\mathcal{S}_d(\R))^{-1}_{\mu_{\beta}^{\otimes d}}$. Moreover, 
	\begin{align*}
		S_{\mu_{\beta}^{\otimes d}}\frac{\text{\rm d}}{\text{\rm d}t}\Phi(t)(\varphi) = \lim_{n\to\infty} S_{\mu_{\beta}^{\otimes d}}\Psi_n(\varphi)=\frac{\text{\rm d}}{\text{\rm d}t}S_{\mu_{\beta}^{\otimes d}}\Phi(t)(\varphi).\tag*{\qed}
	\end{align*}
	
	As the space $\text{\rm Hol}_0\left(\mathcal{S}_d(\R)_\C\right)$ is an algebra, we infer from Equation (\ref{e:st}) that Theorem~\ref{t:chrkd}, Theorem~\ref{t:chrint}, Theorem~\ref{t:chrseq}, and Corollary~\ref{c:dvst} also hold if the $S_{\mu_{\beta}^{\otimes d}}$-transform is replaced by the $T_{\mu_{\beta}^{\otimes d}}$-transform.
	
	\section{Donsker's delta of Mittag-Leffler random vectors}
	In this part, we want to construct a distribution in $(\mathcal{S}_d(\R))^{-1}_{\mu_{\beta}^{\otimes d}}$, $0<\beta<1$, which is a generalization of Donsker's delta of $d$-dimensional Brownian motion. First, for $\eta\in L^2_d(\R)$, define the random vector $G(\cdot,\eta)$ on $\mathcal{S}_d'(\R)$ by
	\begin{align*}
		G(\omega,\eta):= (\langle\omega,\eta_1\mathbf{e}_1\rangle,\dots,\langle\omega,\eta_d\mathbf{e}_d\rangle) =(\langle\omega_1,\eta_1\rangle,\dots,\langle\omega_d,\eta_d\rangle)\in\R^d,\quad\text{for $\mu_{\beta}^{\otimes d}$-a.e.}\ \omega\in\mathcal{S}_d'(\R).
	\end{align*}
	This random vector is well defined as an element of $L^2(\mu_{\beta}^{\otimes d};\R^d)$ by Remark~\ref{r:l2rv}. Let $(\cdot,\cdot)_\text{\rm euc}$ and $|\cdot|_\text{\rm euc}$ be the standard Euclidean scalar product and norm in $\R^d$, respectively. The following properties of $G(\cdot,\eta)$ follow directly from the definition of $\mu_{\beta}^{\otimes d}$ and Proposition~\ref{p:bc01}.
	\begin{prop}\label{p:gprop}
		Let $\eta$, $\zeta\in L^2_d(\R)$ and $p\in\R^d$.
		\begin{enumerate}
			\item[(i)] The characteristic function of $G(\cdot,\eta)$ is given by
			\begin{align*}
				\E_{\mu_{\beta}^{\otimes d}}\left(e^{i(p,G(\cdot,\eta))_{\text{\rm euc}}}\right) = \prod_{k=1}^d E_\beta\left(-\frac{1}{2}p_k^2|\eta_k|_0^2\right).
			\end{align*}
			\item[(ii)] The characteristic function of $G(\cdot,\eta)-G(\cdot,\zeta)$ is given by
			\begin{align*}
				\E_{\mu_{\beta}^{\otimes d}}\left(e^{i(p,G(\cdot,\eta)-G(\cdot,\zeta))_{\text{\rm euc}}}\right) = \prod_{k=1}^d E_\beta\left(-\frac{1}{2}p_k^2|\eta_k-\zeta_k|_0^2\right).
			\end{align*}
			\item[(iii)] The expectation vector of $G(\cdot,\eta)$ is zero, and for all $i,j=1,\dots,d$,
			\begin{align*}
				\E_{\mu_{\beta}^{\otimes d}}\left(G(\cdot,\eta)_i\ G(\cdot,\zeta)_j\right) = \frac{1}{\Gamma(\beta+1)}\delta_{i,j} (\eta_i,\zeta_j)_{L^2(\R)}.
			\end{align*}
			In particular,
			\begin{align*}
				\|\,|G(\cdot,\eta)|_\text{\rm euc}\,\|_{L^2(\mu_{\beta}^{\otimes d})}^2 = \frac{1}{\Gamma(\beta+1)}|\eta|_0^2.
			\end{align*}
			and the covariance matrix of $G(\cdot,\eta)$ is given by $\dfrac{1}{\Gamma(\beta+1)} \text{\rm diag}\left(|\eta_1|_0^2,\dots,|\eta_d|_0^2\right)$.
			\item[(iv)] The components of $G(\cdot,\eta)$ are mutually independent.
		\end{enumerate}
	\end{prop}
	
	\begin{prop}\label{p:srv}
		For $\eta\in L^2_d(\R)$, the $S_{\mu_{\beta}^{\otimes d}}$-transform of the random variable $\langle\cdot,\eta\rangle$ on $\mathcal{S}_d'(\R)$ is given by
		\begin{align}\label{e:srv}
			S_{\mu_{\beta}^{\otimes d}}\langle\cdot,\eta\rangle(\varphi) = \sum_{k=1}^d\frac{E_{\beta,\beta}\left(\frac{1}{2}\langle\varphi_k,\varphi_k\rangle\right)}{\beta E_\beta\left(\frac{1}{2}\langle\varphi_k,\varphi_k\rangle\right)}\langle\varphi_k,\eta_k\rangle
		\end{align}
		for $\varphi\in\mathcal{U}_0$, where $\mathcal{U}_0\subset\mathcal{S}_d(\R)_\C$ is a suitable neighborhood of zero.
	\end{prop}
	\pf Since $\langle\cdot,\eta\rangle\in L^2(\mu_{\beta}^{\otimes d})$ by Remark~\ref{r:l2rv} and Corollary~\ref{c:excom}, for $\varphi:=\varphi^1+i\varphi^2\in\mathcal{U}_0$,
	\begin{align}
		S_{\mu_{\beta}^{\otimes d}}\langle\cdot,\eta\rangle(\varphi) &= \frac{1}{l_{\mu_{\beta}^{\otimes d}}(\varphi)}\int_{\mathcal{S}_d'(\R)}\langle\omega,\eta\rangle e^{\langle\omega,\varphi\rangle}\de\mu_{\beta}^{\otimes d}(\omega)\notag\\
		&= \prod_{k=1}^d \frac{1}{E_{\beta}\left(\frac{1}{2}\langle\varphi_k,\varphi_k\rangle\right)} \int_{\mathcal{S}_d'(\R)} \frac{\text{\rm d}}{\text{\rm d}t}e^{\langle\omega,\varphi\rangle + t\langle\omega,\eta\rangle}\Bigg|_{t=0}\de\mu_{\beta}^{\otimes d}(\omega).\label{e:id}
	\end{align}
	Now, for $t\in[-1,1]$, the map $\mathcal{S}_d'(\R)\ni\omega\mapsto e^{\langle\omega,\varphi\rangle + t\langle\omega,\eta\rangle}$ belongs to $L^1(\mu_{\beta}^{\otimes d})$ by H\"{o}lder's inequality and Corollary~\ref{c:excom}. Moreover, for $\mu_{\beta}^{\otimes d}$-a.e.~$\omega\in\mathcal{S}_d'(\R)$,
	\begin{align*}
		\left|\frac{\text{\rm d}}{\text{\rm d}t}e^{\langle\omega,\varphi\rangle + t\langle\omega,\eta\rangle}\right| = \left|\langle\omega,\eta\rangle\right| |e^{\langle\omega,\varphi\rangle + t\langle\omega,\eta\rangle}| \leq  e^{\langle\omega,\varphi^1\rangle}e^{2|\langle\omega,\eta\rangle|},
	\end{align*}
	and the map $\mathcal{S}_d'(\R)\ni\omega\mapsto e^{\langle\omega,\varphi^1\rangle}e^{2|\langle\omega,\eta\rangle|}$ also belongs to $L^1(\mu_{\beta}^{\otimes d})$ by H\"{o}lder's inequality, Lemma~\ref{p:exp}, and Corollary~\ref{c:excom}. Thus, interchanging the derivative and integral in Equation~(\ref{e:id}) is allowed. Since $e^{\langle\omega,\varphi\rangle + t\langle\omega,\eta\rangle} = e^{\langle\omega,\varphi+t\eta\rangle}$ for each $t\in[-1,1]$ and for $\mu_{\beta}^{\otimes d}$-a.e.~$\omega\in\mathcal{S}_d'(\R)$,  applying Corollary~\ref{c:excom} and Equation~(\ref{e:MLder}) after interchanging the derivative and integral yield
	\begin{align*}
		\int_{\mathcal{S}_d'(\R)} \frac{\text{\rm d}}{\text{\rm d}t}e^{\langle\omega,\varphi\rangle + t\langle\omega,\eta\rangle}\Bigg|_{t=0}\de\mu_{\beta}^{\otimes d}(\omega) &= \frac{\text{\rm d}}{\text{\rm d}t} \left(\prod_{k=1}^d E_\beta\left(\frac{1}{2}\langle\varphi_k+t\eta_k,\varphi_k+t\eta_k\rangle\right)\right)\Bigg|_{t=0}\\
		&= \frac{1}{\beta} \sum_{k=1}^d \langle\varphi_k,\eta_k\rangle E_{\beta,\beta}\left(\frac{1}{2}\langle\varphi_k,\varphi_k\rangle\right)\prod_{j\neq k} E_\beta\left(\frac{1}{2}\langle\varphi_j,\varphi_j\rangle\right).
	\end{align*}
	The last equation and (\ref{e:id}) imply (\ref{e:srv}). \qed
	
	\begin{cor}\label{c:srvc}
		For $\eta\in L^2_d(\R)$ and $\varphi\in\mathcal{U}_0$, $\mathcal{U}_0\subset\mathcal{S}_d(\R)_\C$ a suitable neighborhood of zero,
		\begin{align*}
			S_{\mu_{\beta}^{\otimes d}}G(\cdot,\eta)(\varphi) = \sum_{k=1}^d\frac{E_{\beta,\beta}\left(\frac{1}{2}\langle\varphi_k,\varphi_k\rangle\right)}{\beta E_\beta\left(\frac{1}{2}\langle\varphi_k,\varphi_k\rangle\right)}\langle\varphi_k,\eta_k\rangle \mathbf{e}_k.
		\end{align*}
	\end{cor}
	\pf Apply Proposition~\ref{p:srv} to the random variable $\langle\cdot,\eta_k\mathbf{e}_k\rangle$, $k=1, \dots, d$, and note that $G(\cdot,\eta) = \sum_{k=1}^d\langle\cdot,\eta_k\mathbf{e}_k\rangle\mathbf{e}_k$ for $\mu_{\beta}^{\otimes d}$-almost all $\omega\in\mathcal{S}_d'(\R)$ by Remark~\ref{r:l2rv}. \qed
	
	\begin{thm}\label{t:dkdl}
		Let $0<\beta<1$ and $\eta\in L^2_d(\R)$ such that $\eta_k\neq 0$ for all $k=1,\dots,d$. Then the $d$-dimensional Donsker's delta at $a\in\R^d$, defined via the integral representation
		\begin{align*}
			\delta_a(G(\cdot,\eta)) := \dfrac{1}{(2\pi)^d}\int_{\R^d} e^{i(s,G(\cdot,\eta)-a)_\text{\rm euc}}\de s,
		\end{align*}
		exists in the space $(\mathcal{S}_d(\R))^{-1}_{\mu_{\beta}^{\otimes d}}$ as a weak integral in the sense of Theorem~\ref{t:chrint}.
	\end{thm}
	\pf 
	Since $e^{i(s,G(\cdot,\eta)-a)_\text{\rm euc}}\in L^2(\mu_{\beta}^{\otimes d})$, we apply Corollary~\ref{c:excom} and infer that the map
	\begin{align*}
		\R^d\ni s \mapsto T_{\mu_{\beta}^{\otimes d}}e^{i(s,G(\cdot,\eta)-a)_\text{\rm euc}}(\varphi) = e^{-i(s,a)_\text{\rm euc}}\prod_{k=1}^d E_\beta\left(-\dfrac{1}{2}s_k^2\langle\eta_k,\eta_k\rangle -\dfrac{1}{2}\langle\varphi_k,\varphi_k\rangle - s_k\langle\eta_k,\varphi_k\rangle\right)
	\end{align*}
	is measurable for all $\varphi\in\mathcal{S}_d(\R)_\C$. Now, by Proposition 5.2 in \cite{GJRS15}, for each $k=1,\dots,d$, there exists a constant $C_k<\infty$ such that 
	$$
		\int_{\R} \left|E_\beta\left(-\dfrac{1}{2}s_k^2\langle\eta_k,\eta_k\rangle -\dfrac{1}{2}\langle\phi,\phi\rangle - s_k\langle\eta_k,\phi\rangle\right)\right|\de s_k \leq C_k
	$$
	for all $\phi$ belonging to the set $\{\phi\in \mathcal{S}(\R)_\C: |\phi|_0<M \}$ for any $0<M<\infty$. Hence, for all $\varphi$ belonging to $\mathcal{U}_0:=\{\varphi\in \mathcal{S}_d(\R)_\C: |\varphi|_0<M \}$, $0<M<\infty$,
	\begin{align*}
		\dfrac{1}{(2\pi)^d}\int_{\R^d} \big|T_{\mu_{\beta}^{\otimes d}}e^{i(s,G(\cdot,\eta)-a)_\text{\rm euc}}(\varphi)\big|\de s &\leq \dfrac{1}{(2\pi)^d} \prod_{k=1}^d C_k <\infty.
	\end{align*}
	Therefore, $\delta_a(\langle\cdot,\eta\rangle)\in (\mathcal{S}_d(\R))^{-1}_{\mu_{\beta}^{\otimes d}}$ by Theorem~\ref{t:chrint}. \qed
	
	\begin{rmk}\label{r:tdd}
		We can use Theorem 5.3 in \cite{GJRS15} to obtain an explicit formula for the $T_{\mu_{\beta}^{\otimes d}}$-transform of the Donsker's delta at $0$: for all $\varphi\in\mathcal{U}_0$, $\mathcal{U}_0$ as in the proof of Theorem~\ref{t:dkdl},
		\begin{align*}
			T_{\mu_{\beta}^{\otimes d}}\delta_0(G(\cdot,\eta))(\varphi) = \dfrac{1}{(2\pi)^{d/2}}\prod_{k=1}^d\langle\eta_k,\eta_k\rangle^{-1/2} H^{1\,1}_{1\,2}\left(\frac{1}{2}\langle\varphi_k,\varphi_k\rangle-\dfrac{\langle\eta_k,\varphi_k\rangle}{\langle\eta_k,\eta_k\rangle}\Bigg| \begin{matrix} (\tfrac{1}{2},1)\\ (0,1),(\tfrac{1}{2}\beta,\beta) \end{matrix} \right),
		\end{align*}
		where $H$ is the Fox $H$-function (see~Appendix A in \cite{GJRS15}). In particular, the generalized expectation of the Donsker's delta at $0$ is given by
		\begin{align*}
			\E_{\mu_{\beta}^{\otimes d}}(\delta_0(G(\cdot,\eta))) := \llangle \delta_0(G(\cdot,\eta)),1 \rrangle_{\mu_{\beta}^{\otimes d}} = T_{\mu_{\beta}^{\otimes d}}\delta_0(G(\cdot,\eta))(0) = \frac{1}{2^{d/2}\Gamma(1-\tfrac{1}{2}\beta)^d} \prod_{k=1}^d\langle\eta_k,\eta_k\rangle^{-1/2}.
		\end{align*}
	\end{rmk}
	
	\section{Vector-valued generalized grey Brownian motion}
	
	\subsection{Real-valued generalized grey Brownian motion}
	Here, we discuss the definition of the generalized grey Brownian motion introduced by Schneider in \cite{Sch92}. We follow the construction from \cite{GJ16}, starting with the Mittag-Leffler measure $\mu_\beta$ on $\mathcal{S}'(\mathbb{R})$. Let $0<\alpha<2$ be given and define the operator $M^{\alpha/2}_\pm$ on $\mathcal{S}(\mathbb{R})$ as follows:
	\begin{align*}
		\mathcal{S}(\mathbb{R})\ni \varphi \mapsto M^{\alpha/2}_\pm\varphi := 
		\left\{
		\begin{array}{cl}
			K_{\alpha/2} D^{(1-\alpha)/2}_\pm\varphi, & 0<\alpha<1,\\[5pt]
			\varphi, & \alpha = 1,\\[5pt]
			K_{\alpha/2} I^{(\alpha-1)/2}_\pm\varphi, & 1<\alpha<2,
		\end{array}
		\right.
	\end{align*}
	where $K_{\alpha/2}:=\sqrt{\alpha\sin(\alpha\pi/2)\Gamma(\alpha)}$ is a normalization constant, and for $r>0$, $D^r_\pm$ is the 
	right-sided and the left-sided Marchaud fractional derivative of order $r$, while $I^r_\pm$ denote the right-sided and the left-sided  Riemann-Liouville fractional integral of order $r$. See \cite{Ben03,SKM93} for further details of these operators. Although $M^{\alpha/2}_\pm$ is defined on $\mathcal{S}(\mathbb{R})$, its domain is larger. In particular, the domain includes the indicator function $\mathbbm{1}_{[0,t)}$, $t\geq 0$, and that $M^{\alpha/2}_\pm\mathbbm{1}_{[0,t)}\in L^2(\mathbb{R})$ (see Remark~3.2 in \cite{GJ16}). Moreover, the following scalar product holds.
	\begin{prop}\label{p:malpf}{\rm\bf(Corollary~3.5 in \cite{GJ16})}
		For all $t,s\geq 0$, $\alpha\in(0,2)$ and $m,n\in\mathbb{N}$,
		\begin{align}\label{e:morth}
			(M^{\alpha/2}_{-}\mathbbm{1}_{[0,t)},M^{\alpha/2}_{-}\mathbbm{1}_{[0,s)})_{L^2(\R)} = \frac{1}{2}\left(t^\alpha+s^\alpha-|t-s|^\alpha\right).
		\end{align}
	\end{prop}
	
	\begin{defn} For $\beta\in(0,1]$, $\alpha\in(0,2)$ and $t\geq 0$, define $B^{\beta,\alpha}_t$ as follows:
		\begin{align*}
			\mathcal{S}'(\mathbb{R})\ni \omega \mapsto B^{\beta,\alpha}_t(\omega):=\langle\omega,M^{\alpha/2}_{-}\mathbbm{1}_{[0,t)}\rangle,
		\end{align*}
		The process $B^{\beta,\alpha}:=(B^{\beta,\alpha}_t)_{t\geq 0}$ takes values in $L^2(\mu_\beta)$, and is called a \textit{generalized grey Brownian motion} (briefly \textit{ggBm}). If $\alpha=\beta$, the process $B^{\beta,\beta}$ is denoted by $B^{\beta}:=(B^{\beta}_t)_{t\geq 0}$, and is called a \textit{grey Brownian motion} (briefly \textit{gBm}). 
	\end{defn}

	\begin{rmk}
		The family of processes $B^{\beta,\alpha}$, $\beta\in(0,1]$, $\alpha\in(0,2)$, includes the following:
		\begin{itemize}
			\item The process $B^{1}=B^{1,1}$ is a standard one-dimensional Brownian motion (see \cite{HKPS93}).
			\item The process $B^{1,\alpha}$ is a one-dimensional fractional Brownian motion with Hurst parameter $\alpha/2$ (see \cite{Ben03}).
		\end{itemize}
	\end{rmk}

	We state the following properties of the process $B^{\beta,\alpha}$ from \cite{GJ16}.
	
	\begin{prop} Let $\beta\in(0,1]$ and $\alpha\in(0,2)$.
		\begin{enumerate}
			\item[(i)] $B^{\beta,\alpha}$ has mean zero and covariance
			\begin{align*}
				\mathbb{E}_{\mu_\beta}\left(B^{\beta,\alpha}_t B^{\beta,\alpha}_s\right) = \dfrac{1}{2\,\Gamma(\beta+1)}\left(t^\alpha+s^\alpha-|t-s|^\alpha\right),\quad t,s\geq 0.
			\end{align*}
			\item[(ii)] For all $p\in\mathbb{N}$, there exists $K<\infty$ such that
			\begin{align*}
				\mathbb{E}_{\mu_\beta}\left(\left|B^{\beta,\alpha}_t - B^{\beta,\alpha}_s\right|^{2p}\right) \leq K|t-s|^{\alpha p},\quad t,s\geq 0,
			\end{align*}
			and hence $B^{\beta,\alpha}$ has a continuous modification by Kolmogorov's continuity theorem.
			\item[(iii)] $B^{\beta,\alpha}$ has stationary increments.
		\end{enumerate}
	\end{prop}
	
	\subsection{Construction, properties and its noise}
	
	\begin{defn}
		For $d\in\N$, $0<\beta\leq 1$, $0<\alpha<2$, and $t\geq 0$, define $B^{\beta,\alpha}_{d,t}$ as follows:
		\begin{align*}
			\mathcal{S}_d'(\R) \ni \omega \mapsto B^{\beta,\alpha}_{d,t}(\omega) := (\langle\omega_1,M^{\alpha/2}_{-}\mathbbm{1}_{[0,t)}\rangle,\dots,\langle\omega_d,M^{\alpha/2}_{-}\mathbbm{1}_{[0,t)}\rangle)\in\R^d.
		\end{align*}
		The process $B^{\beta,\alpha}_{d}:=(B^{\beta,\alpha}_{d,t})_{t\geq 0}$ takes values in $L^2(\mu_{\beta}^{\otimes d};\R^d)$, and is called a \textit{vector-valued generalized grey Brownian motion} (briefly \textit{vggBm}). If $\alpha=\beta$, the process $B^{\beta,\beta}_{d}$ is denoted by $B^{\beta}_{d}:=(B^{\beta}_{d,t})_{t\geq 0}$, and is called a \textit{vector-valued grey Brownian motion} (briefly \textit{vgBm}).
	\end{defn}
	
	If we define a vggBm using the standard Mittag-Leffler measure $\mu_{\beta}$ on $\mathcal{S}_d'(\R)$, $d\geq 2$, then in view of Proposition~\ref{p:indmub}, it has independent components if and only if it is a fractional Brownian motion. However, since this process is defined using the product measure $\mu_{\beta}^{\otimes d}$, we get the following result.
	\begin{prop}
		The process $B^{\beta,\alpha}_{d}$ has the following properties.
		\begin{enumerate}
			\item[(i)] For each $t\geq 0$, $B^{\beta,\alpha}_{d,t}$ has characteristic function
			\begin{align*}
				\E_{\mu_{\beta}^{\otimes d}}\left(e^{i(p,B^{\beta,\alpha}_{d,t})_{\text{\rm euc}}}\right) = \prod_{k=1}^d E_\beta\left(-\frac{|p|_\text{\rm euc}^2}{2}t^\alpha\right).
			\end{align*}
			\item[(ii)] For each $t\geq 0$, $B^{\beta,\alpha}_{d,t}$ has expectation zero, and for all $i,j=1,\dots,d$ and $t,s\geq 0$,
			\begin{align*}
				\E_{\mu_{\beta}^{\otimes d}}\left((B^{\beta,\alpha}_{d,t})_i(B^{\beta,\alpha}_{d,s})_j\right) = \dfrac{1}{2\,\Gamma(\beta+1)}\delta_{i,j}\left(t^\alpha+s^\alpha-|t-s|^\alpha\right).
			\end{align*}
			In particular, for each $t\geq 0$, the covariance matrix of $B^{\beta,\alpha}_{d,t}$ is given by $\dfrac{t^\alpha}{\Gamma(\beta+1)} I_d$, where $I_d$ is the identity matrix of order $d$.
			\item[(iii)] For each $t\geq 0$, $B^{\beta,\alpha}_{d,t}$ has independent components.
			\item[(iv)] $B^{\beta,\alpha}_{d}$ is $\alpha/2$ self-similar with stationary increments in the strict sense.
		\end{enumerate}
	\end{prop}
	\pf Statements (i)-(iii) follow directly from Proposition~\ref{p:gprop} and Equation~\ref{e:morth}. Let $p:=(p_{k,r})\in\mathbb{R}^{d\times n}$ and $0\leq t_{1}<t_{2}<\ldots<t_{n}$. To show $\alpha/2$ self-similarity, we need to prove that for all $a>0$,
	\begin{align}
		\E_{\mu_{\beta}^{\otimes d}}\left(\exp\left(i\sum_{r=1}^{n}\sum_{k=1}^{d}p_{k,r}\left\langle (\cdot)_k,M_{-}^{\alpha/2}\mathbbm{1}_{[0,at_{r})}\right\rangle \right)\right) &= \E_{\mu_{\beta}^{\otimes d}}\left(\exp\left(ia^{\alpha/2}\sum_{r=1}^{n}\sum_{k=1}^{d}p_{k,r}\left\langle (\cdot)_k,M_{-}^{\alpha/2}\mathbbm{1}_{[0,t_{r})}\right\rangle \right)\right). \label{e:ecfss}
	\end{align}
	By Proposition~\ref{p:gprop}(i), Equation (\ref{e:ecfss}) is equivalent to
	\begin{align*}
		\prod_{k=1}^{d} E_{\beta}\left(-\frac{1}{2}\left|\sum_{r=1}^{n}p_{k,r}M_{-}^{\alpha/2}\mathbbm{1}_{[0,at_{r})}\right|_{0}^{2}\right) = \prod_{k=1}^{d} E_{\beta}\left(-\frac{1}{2}\left|a^{\alpha/2}\sum_{r=1}^{n}p_{k,r}M_{-}^{\alpha/2}\mathbbm{1}_{[0,t_{r})}\right|_{0}^{2}\right),
	\end{align*}
	and this equation holds. Indeed, for all $k=1,\dots,d$, we can use Equation (\ref{e:morth}) to infer that
	\[
	\left|\sum_{r=1}^{n}p_{k,r}M_{-}^{\alpha/2}\mathbbm{1}_{[0,at_{r})}\right|_{0}^{2} =  a^{\alpha}\left|\sum_{r=1}^{n}p_{k,r}M_{-}^{\alpha/2}\mathbbm{1}_{[0,t_{r})}\right|_{0}^{2}.
	\]
	A similar procedure may be applied in order to prove that the increments are stationary; we have to show that for all $h\geq 0$,
	\begin{multline}\label{e:ecfst}
		\E_{\mu_{\beta}^{\otimes d}}\left(\exp\left(i\sum_{r=1}^{n}\sum_{k=1}^{d}p_{k,r}\left\langle (\cdot)_k,M_{-}^{\alpha/2}\mathbbm{1}_{[0,t_{r}+h)} - M_{-}^{\alpha/2}\mathbbm{1}_{[0,h)}\right\rangle \right)\right) \\= 
		\E_{\mu_{\beta}^{\otimes d}}\left(\exp\left(i\sum_{r=1}^{n}\sum_{k=1}^{d}p_{k,r}\left\langle (\cdot)_k,M_{-}^{\alpha/2}\mathbbm{1}_{[0,t_{r})}\right\rangle \right)\right).
	\end{multline}
	Equation (\ref{e:ecfst}) holds if for all $k=1,\dots,d$,
	\[
	\left|\sum_{r=1}^{n}p_{k,r}(M_{-}^{\alpha/2}\mathbbm{1}_{[0,t_{r}+h)} - M_{-}^{\alpha/2}\mathbbm{1}_{[0,h)})\right|_{0}^{2} = \left|\sum_{r=1}^{n}p_{k,r}M_{-}^{\alpha/2}\mathbbm{1}_{[0,t_{r})}\right|_{0}^{2},
	\]
	which can be verified analogously as that of self-similarity. \qed
	\bigskip
	
	To determine the derivative of the vggBm $B^{\beta,\alpha}_{d}$ in the sense of Corollary~\ref{c:dvst}, we use Corollary~\ref{c:srvc} and Equation~(21) of \cite{GJ16} to infer that for all $t\geq 0$,
	\begin{align*}
		S_{\mu_{\beta}^{\otimes d}}B^{\beta,\alpha}_{d,t}(\varphi) &= \sum_{k=1}^d\frac{E_{\beta,\beta}\left(\frac{1}{2}\langle\varphi_k,\varphi_k\rangle\right)}{\beta E_\beta\left(\frac{1}{2}\langle\varphi_k,\varphi_k\rangle\right)}\langle\varphi_k,M^{\alpha/2}_{-}\mathbbm{1}_{[0,t)}\rangle \mathbf{e}_k\\
		&= \sum_{k=1}^d\frac{E_{\beta,\beta}\left(\frac{1}{2}\langle\varphi_k,\varphi_k\rangle\right)}{\beta E_\beta\left(\frac{1}{2}\langle\varphi_k,\varphi_k\rangle\right)}\int_{0}^{t}\left(M^{\alpha/2}_{+}\varphi_k\right)(x)\de x\ \mathbf{e}_k,
	\end{align*}
	on the set $\mathcal{U}_\beta:=\{\varphi\in\mathcal{S}_d(\R)_\C:|G(\varphi,\varphi)|_{\text{euc}}<\varepsilon_\beta\}$, where $\varepsilon_\beta>0$ is chosen so that $E_\beta(z)>0$ for all $|z|<\varepsilon_\beta$. We infer from the continuity of $M^{\alpha/2}_{+}\varphi_k$ on $\R$ (see Theorem~2.7 in \cite{Ben03}) that for $\varphi\in\mathcal{U}_\beta$,
	\begin{align*}
		\frac{\text{\rm d}}{\text{\rm d}t}S_{\mu_{\beta}^{\otimes d}}B^{\beta,\alpha}_{d,t}(\varphi) = \sum_{k=1}^d\frac{E_{\beta,\beta}\left(\frac{1}{2}\langle\varphi_k,\varphi_k\rangle\right)}{\beta E_\beta\left(\frac{1}{2}\langle\varphi_k,\varphi_k\rangle\right)}\big(M^{\alpha/2}_{+}\varphi_k\big)(t)\ \mathbf{e}_k.
	\end{align*}
	Now, as the Mittag-Leffler functions are holomorphic, there are $p,q\in\N$ and a constant $K<\infty$ such that $\mathcal{U}_{p,q}\subset\mathcal{U}_\beta$, and for all $\varphi\in\mathcal{U}_{p,q}$, and $k=1,\dots,d$,
	\begin{align*}
		\left|\frac{E_{\beta,\beta}\left(\frac{1}{2}\langle\varphi_k,\varphi_k\rangle\right)}{\beta E_\beta\left(\frac{1}{2}\langle\varphi_k,\varphi_k\rangle\right)}\right| \leq K.
	\end{align*}
	Moreover, Theorem~2.3 in \cite{Ben03} shows that there exists $p'\in\N$ and a constant $C<\infty$ such that for all $k=1,\dots,d$, and $x\in\R$, 
	\begin{align*}
		\left|\big(M^{\alpha/2}_{+}\varphi_k\big)(x)\right| \leq C|\varphi_k|_{p'}.
	\end{align*}
	Thus, by choosing $p^*>\max\{p,p'\}$, the following estimate holds for all $t\geq 0$, $\varphi\in\mathcal{U}_{p^*,q}$, and $k=1,\dots,d$,
	\begin{align}\label{e:CBound}
		\left|\frac{E_{\beta,\beta}\left(\frac{1}{2}\langle\varphi_k,\varphi_k\rangle\right)}{\beta E_\beta\left(\frac{1}{2}\langle\varphi_k,\varphi_k\rangle\right)}\big(M^{\alpha/2}_{+}\varphi_k\big)(t)\right| \leq KC|\varphi|_{p^*} \leq 2^{-q}KC.
	\end{align}
	Therefore, by Corollary~\ref{c:dvst}, we establish the existence of the derivative of $B^{\beta,\alpha}_{d,t}$.
	\begin{prop}\label{p:dervggBm}
		For each $t\geq 0$, $\frac{\text{\rm d}}{\text{\rm d}t}B^{\beta,\alpha}_{d,t}$ exists as a vector with components in $(\mathcal{S}_d(\R))^{-1}_{\mu_{\beta}^{\otimes d}}$ in the sense of Corollary~\ref{c:dvst}. Moreover, for all $\varphi$ belonging to a suitable neighborhood of zero in $\mathcal{S}_d(\R)_\C$,
		\begin{align*}
			S_{\mu_{\beta}^{\otimes d}}\frac{\text{\rm d}}{\text{\rm d}t}B^{\beta,\alpha}_{d,t}(\varphi) = \sum_{k=1}^d\frac{E_{\beta,\beta}\left(\frac{1}{2}\langle\varphi_k,\varphi_k\rangle\right)}{\beta E_\beta\left(\frac{1}{2}\langle\varphi_k,\varphi_k\rangle\right)}\big(M^{\alpha/2}_{+}\varphi_k\big)(t)\ \mathbf{e}_k,\quad t\geq 0.
		\end{align*}
	\end{prop}	
	\subsection{Local time and self-intersection local time for vggBm}
	Here, we consider the local time and the self-intersection local time for vggBm, which are given respectively by
		\begin{gather*}
			L_{\beta,\alpha}(a,T) := \int_{[0,T]} \delta_a(B^{\beta,\alpha}_{d,t})\de t,\quad a\in\R^d,\ T>0,\\
			L_{\beta,\alpha}^s(T) :=\int_{[0,T]}\int_{[0,T]}\delta_0(B^{\beta,\alpha}_{d,s}-B^{\beta,\alpha}_{d,u})\de u\de s,\quad T>0.
		\end{gather*}
		The vggBm local time $L_{\beta,\alpha}(a,T)$ is used to measure the amount of time the sample path of a vggBm spends at $a\in\R^d$ within the time interval $[0,T]$, while the vggBm self-intersection local time $L_{\beta,\alpha}^s(T)$ is intended to measure the amount of time in which the sample path of a vggBm spends intersecting itself also within the time interval $[0,T]$. A priori, the expressions above have no mathematical meaning, since Lebesgue integration of Dirac delta distribution is not defined. In the following, we prove that under some constraints, we can make sense of these objects as weak integrals in the sense of Theorem~\ref{t:chrint}.
	\begin{thm}\label{p:ltvggBm}
		For $d\in\N$, $0<\beta<1$, $0<\alpha<2/d$, $T>0$ and $a\in\R^d$, the vggBm local time $L_{\beta,\alpha}(a,T)$
		and the vggBm self-intersection local time $L_{\beta,\alpha}^s(T)$ exist in $(\mathcal{S}_d(\R))^{-1}_{\mu_{\beta}^{\otimes d}}$ as weak integrals in the sense of Theorem~\ref{t:chrint}. Moreover,
		for all $\varphi$ belonging to a suitable neighborhood $\mathcal{U}_0\subset\mathcal{S}_d(\R)_\C$ of zero,
		\begin{gather*}
			T_{\mu_{\beta}^{\otimes d}}L_{\beta,\alpha}(a,T)(\varphi) = \int_{[0,T]}T_{\mu_{\beta}^{\otimes d}} \delta_a(B^{\beta,\alpha}_{d,t})(\varphi)\de t.\\
			T_{\mu_{\beta}^{\otimes d}}L_{\beta,\alpha}^s(T)(\varphi) = \int_{[0,T]}\int_{[0,T]}T_{\mu_{\beta}^{\otimes d}}\delta_0(B^{\beta,\alpha}_{d,s}-B^{\beta,\alpha}_{d,u})(\varphi)\de u\de s.
		\end{gather*}
	\end{thm}
	\pf Let $\varphi\in\mathcal{S}_d(\R)_\C$ with $|\varphi|<M$, for some $M<\infty$. Following the same calculations from Proposition~5.2 in \cite{GJRS15} and using Proposition~\ref{p:malpf}, we have
	\begin{align*}
		&\int_{[0,T]}\left|T_{\mu_{\beta}^{\otimes d}}\delta_a(B^{\beta,\alpha}_{d,t})(\varphi)\right|\de t\\
		&\hspace{1cm} \leq \frac{1}{(2\pi)^d} \int_{[0,T]} \prod_{k=1}^d\int_{\R}\left|E_\beta\left(-\dfrac{1}{2}s^2\langle M^{\alpha/2}_{-}\mathbbm{1}_{[0,t)}, M^{\alpha/2}_{-}\mathbbm{1}_{[0,t)}\rangle -\dfrac{1}{2}\langle\varphi_k,\varphi_k\rangle - s\langle M^{\alpha/2}_{-}\mathbbm{1}_{[0,t)},\varphi_k\rangle\right)\right|\de s\de t\\
		&\hspace{1cm} \leq \frac{1}{(2\pi)^d} \int_{[0,T]} \prod_{k=1}^d \left[\sqrt{\frac{2\pi}{t^\alpha}}\int_{0}^{\infty} M_\beta(r)r^{-1/2}\exp\left(\frac{1}{2}M^2r\right)\de r\right]\de t,
	\end{align*}
	where $M_\beta$ is the $M$-Wright function (see~\cite{Mai95}). By Lemma~A.4 in \cite{GJRS15},
	$$
		K:=\int_{0}^{\infty} M_\beta(r)r^{-1/2}\exp\left(\frac{1}{2}M^2r\right)\de r <\infty,
	$$
	so that
	\begin{align*}
		\int_{[0,T]}\left|T_{\mu_{\beta}^{\otimes d}}\delta_a(B^{\beta,\alpha}_{d,t})(\varphi)\right|\de t \leq \frac{K^d}{(2\pi)^{d/2}}
		\int_{[0,T]} t^{-\alpha d/2}\de t = \frac{2}{2-\alpha d}K^dT^{1-\alpha d/2} <\infty.
	\end{align*}
	Therefore, $L_{\beta,\alpha}(a,T)\in (\mathcal{S}_d(\R))^{-1}_{\mu_{\beta}^{\otimes d}}$ by Theorem~\ref{t:chrint}. A similar computation holds for the case of $L_{\beta,\alpha}^s(T)$: if we set $\eta_{x}:=M_{-}^{\alpha/2}\mathbbm{1}_{[0,x)}$ for $x\in[0,T]$, then for all $\varphi\in\mathcal{S}_d(\mathbb{R})_\mathbb{C}$ with $|\varphi|<M$, $M<\infty$, and $s,u\in[0,T]$,
	\begin{align*}
		\left|T_{\mu_{\beta}^{\otimes d}}\delta(B^{\beta,\alpha}_{d,s}-B^{\beta,\alpha}_{d,u})(\varphi)\right|
		&\leq \frac{1}{(2\pi)^d} \prod_{k=1}^d\int_{\mathbb{R}}\left|E_\beta\left(-\dfrac{1}{2}\lambda_k^2\langle \eta_s-\eta_u,\eta_s-\eta_u\rangle -\dfrac{1}{2}\langle\varphi_k,\varphi_k\rangle - \lambda_k\langle \eta_s-\eta_u,\varphi_k\rangle\right)\right|\de\lambda_k\\
		&\leq \frac{1}{(2\pi)^d} \prod_{k=1}^d\left(\sqrt{\frac{2\pi}{\langle \eta_s-\eta_u,\eta_s-\eta_u\rangle}}\int_{0}^{\infty} M_\beta(r)r^{-1/2}\exp\left(\frac{1}{2}M^2r\right)\de r\right)\\
		&= \frac{K^d}{(2\pi)^{d/2}}|s-u|^{-\alpha d/2},
	\end{align*}
	so that
	\begin{align*}
		\int_{[0,T]}\int_{[0,T]}\left|T_{\mu_{\beta}^{\otimes d}}\delta(B^{\beta,\alpha}_{d,s}-B^{\beta,\alpha}_{d,u})(\varphi)\right|\de u\de s &\leq \frac{2K^d}{(2\pi)^{d/2}}\int_{[0,T]}\int_{[0,s]}(s-u)^{-\alpha d/2}\de u\de s\\
		&= \frac{8K^dT^{2-\alpha d/2}}{(2\pi)^{d/2}(2-\alpha d)(4-\alpha d)} <\infty.
	\end{align*}
	The conclusion follows from Theorem~\ref{t:chrint}. \qed

	\begin{rmk} Let $d\in\N$ and $0<\alpha<2/d$.
		\begin{enumerate}
			\item[(i)] If $0<\beta<1$, then by Theorem~\ref{p:ltvggBm} and Remark~\ref{r:tdd}, the generalized expectation of $L_{\beta,\alpha}(0,T)$ is given by
			\begin{align}\label{e:explt}
				\E_{\mu_{\beta}^{\otimes d}}(L_{\beta,\alpha}(0,T)) = \int_{[0,T]} T_{\mu_{\beta}^{\otimes d}}\delta_0(B^{\beta,\alpha}_{d,t})(0)\de t = \frac{T^{1-\alpha d/2}}{2^{d/2-1}\Gamma(1-\tfrac{1}{2}\beta)^d(2-\alpha d)},
			\end{align}
			while the generalized expectation of $L_{\beta,\alpha}^s(T)$ is given by
			\begin{align}\label{e:expslt}
				\E_{\mu_{\beta}^{\otimes d}}(L_{\beta,\alpha}^s(T)) = \int_{[0,T]}\int_{[0,T]}T_{\mu_{\beta}^{\otimes d}}\delta_0(B^{\beta,\alpha}_{d,s}-B^{\beta,\alpha}_{d,u})(0)\de u\de s = \frac{T^{2-\alpha d/2}}{2^{d/2-2}\Gamma(1-\tfrac{1}{2}\beta)^d(2-\alpha d)(4-\alpha d)}.
			\end{align}
			\item[(ii)] Consider the case $\beta=1$, in which the process $B^{\beta,\alpha}_d$ is a $d$-dimensional fractional Brownian motion (fBm) with Hurst parameter $H=\alpha/2$. In this case, the assumption that $\alpha d<2$ reduces to $Hd<1$, and Corollary~4.10(a) in \cite{HO02} shows that the right-hand side of (\ref{e:explt}) corresponds to the generalized expectation of the fBm local time at $0$. Moreover, a simple application of Lebesgue's dominated convergence theorem to Equation (14) in \cite{HN05} shows that the right-hand side of (\ref{e:expslt}) corresponds to the expectation of the $L^2$-limit of the approximated self-intersection local time $I_\varepsilon$ of fBm defined by Equation~(2) in \cite{HN05}.
		\end{enumerate}
	\end{rmk}
	
	\subsection{Linear stochastic differential equations driven by vggBm noise}
	In this section, we study linear stochastic differential systems of the form 
	\begin{equation}\label{e:SDE}
		\left\{
		\begin{aligned}
			\text{\rm d}X_t &= A(t)X_t\,{\text{\rm d}t} + \sigma\,\text{\rm d}{B_{d,t}^{\beta, \alpha}}, \quad t\in [0,T]\\
			X_0 &=x_0 \in \mathbb{R}^d,
		\end{aligned}
		\right.
	\end{equation}
	where we assume that for $T>0$, $A:[0,T]\to\R^{d\times d}$ is continuous, and $\sigma\in\R$. As in the case for white noise analysis, we rewrite (\ref{e:SDE}) as a system of equations in $(\mathcal{S}_d(\R))^{-1}_{\mu_{\beta}^{\otimes d}}$:
	\begin{align}\label{e:DSDE}
		\left\{
		\begin{aligned}
			\frac{\text{\rm d}}{\text{\rm d}t} X_t &= A(t)  X_t + \sigma\,\frac{\text{\rm d}}{\text{\rm d}t}{B_{d,t}^{\beta, \alpha}}, \quad t\in [0,T]\\
			X_0 &=x_0 \in \mathbb{R}^d,
		\end{aligned}
		\right.
	\end{align}
	and seek a vector-valued process $X_t$ with components taking values in $(\mathcal{S}_d(\R))^{-1}_{\mu_{\beta}^{\otimes d}}$ that solves system (\ref{e:DSDE})  for all $t\in [0,T]$.
	
	First, assume that there exists such a process $X_t$. If we apply the $S_{\mu_{\beta}^{\otimes d}}$-transform to both sides of (\ref{e:DSDE})$_1$, then by Corollary~\ref{c:dvst}, for some neighborhood $\mathcal{U}_0$ of zero in $\mathcal{S}_d(\R)_\C$,
	\begin{align}\label{e:SDSDE}
		\frac{\text{\rm d}}{\text{\rm d}t} S_{\mu_{\beta}^{\otimes d}}X_t(\varphi) = A(t) S_{\mu_{\beta}^{\otimes d}}X_t(\varphi) + \sigma S_{\mu_{\beta}^{\otimes d}}\frac{\text{\rm d}}{\text{\rm d}t} B_{d,t}^{\beta,\alpha}(\varphi), 
		\quad t\in [0,T],\ \varphi \in \mathcal{U}_0.
	\end{align}
	Note that in (\ref{e:SDSDE}), the matrix $A(t)$ and $S_{\mu_{\beta}^{\otimes d}}$-transform commute since $A(t)$ is independent of $\omega\in\mathcal{S}_d'(\R)$. Set $Y_t(\varphi):=S_{\mu_{\beta}^{\otimes d}}X_t(\varphi)$ and use Proposition~\ref{p:dervggBm} to obtain
	\begin{align}\label{e:ODE}
		\frac{\text{\rm d}}{\text{\rm d}t} Y_t(\varphi) = A(t) Y_t(\varphi) +   \sigma C_{\beta,\alpha}(\varphi,t) , \quad t\in [0,T],\ \varphi \in \mathcal{U}_0,
	\end{align}
	where, for convenience, we set
	$$
		C_{\beta,\alpha}(\varphi,t) := \sum_{k=1}^d\frac{E_{\beta,\beta}\left(\frac{1}{2}\langle\varphi_k,\varphi_k\rangle\right)}{\beta E_\beta\left(\frac{1}{2}\langle\varphi_k,\varphi_k\rangle\right)}\big(M^{\alpha/2}_{+}\varphi_k\big)(t)\,\mathbf{e}_k.
	$$
	Equation~(\ref{e:ODE}) is a linear nonhomogeneous ordinary differential system, with initial condition
	\begin{align}\label{e:ODEin}
		Y_0(\varphi) = S_{\mu_{\beta}^{\otimes d}}X_0(\varphi) = x_0.
	\end{align}
	It has a unique solution on $[0,T]$ for each $\varphi \in \mathcal{U}_0$, since both $A$ and $M^{\alpha/2}_{+}\varphi_k$, $k=1,\dots,d$, are continuous on $[0,T]$. The solution of (\ref{e:ODE}) with initial condition (\ref{e:ODEin}) is computed using the method of variation of constants:
	\begin{align}\label{e:SODE}
		Y_t(\varphi) &=V(t)V(0)^{-1}x_0 + \sigma V(t)\int_{0}^{t} V(s)^{-1} C_{\beta,\alpha}(\varphi,s)\de s, \quad t\in [0,T],\ \varphi \in \mathcal{U}_0,
	\end{align}
	where $V:[0,T]\to\R^{d\times d}$ is a fundamental matrix to the homogeneous system
	\begin{align}\label{e:hODE}
		\frac{\text{\rm d}}{\text{\rm d}t} \mathbf{y}(t) = A(t) \mathbf{y}(t),\quad t\in[0,T].
	\end{align}

	Now, let $u_{j,k}:\R\to\R$, $j,k=1,\dots,d$, be the $(j,k)$-entry of $V^{-1}$, extended to zero outside $[0,T]$. Since $u_{j,k}$ is continuously differentiable on the compact interval $[0,T]$, for each $t\in[0,T]$, the product $\mathbbm{1}_{[0,t)}u_{j,k}$ belongs to $L^q(\R)$ for every $1\leq q\leq \infty$. Moreover, $M^{\alpha/2}_{-}\left(\mathbbm{1}_{[0,t)}u_{j,k}\right)\in L^2(\R)$. Indeed, this is clear for $\alpha=1$. If $1<\alpha<2$, then the statement follows from Theorem~5.3 in \cite{SKM93}, since $\mathbbm{1}_{[0,t)}u_{j,k}\in L^{2/\alpha}(\R)$. For $0<\alpha<1$, the function $\mathbbm{1}_{[0,t)}u_{j,k}$ is piecewise Lipschitz with a finite number of discontinuities and $\mathbbm{1}_{[0,t)}(x)u_{j,k}(x)\to 0$ as $|x|\to\infty$, so that $M^{\alpha/2}_{-}\left(\mathbbm{1}_{[0,t)}u_{j,k}\right)\in L^2(\R)$ by Theorem~11.7 and Theorem 6.1 in \cite{SKM93}. Furthermore, by following a proof similar to that of Lemma 2.5 in \cite{Ben03}, we have the following duality relation:
	\begin{align}\label{e:mdual}
		\langle\phi,M^{\alpha/2}_{-}\left(\mathbbm{1}_{[0,t)}u_{j,k}\right)\rangle = \langle M^{\alpha/2}_{+}\phi,\mathbbm{1}_{[0,t)}u_{j,k}\rangle = \int_{0}^{t} (M^{\alpha/2}_{+}\phi)(s)\,u_{j,k}(s)\de s, \quad \phi\in\mathcal{S}(\R)_\C.
	\end{align}
	For each $t\in[0,T]$ and $j=1,\dots,d$, set
	$$
		\eta_{t,j} := \sum_{k=1}^{d}M^{\alpha/2}_{-}\left(\mathbbm{1}_{[0,t)}u_{j,k}\right)\mathbf{e}_k.
	$$
	Then for all $t\in[0,T]$ and $j=1,\dots,d$, $\eta_{t,j}$ belongs to $L^2_d(\R)$. 
	Moreover, by Proposition~\ref{p:srv} and Equation (\ref{e:mdual}),
	\begin{align*}
		S_{\mu_{\beta}^{\otimes d}}\langle\cdot,\eta_{t,j}\rangle(\varphi) &= \sum_{k=1}^d\frac{E_{\beta,\beta}\left(\frac{1}{2}\langle\varphi_k,\varphi_k\rangle\right)}{\beta E_\beta\left(\frac{1}{2}\langle\varphi_k,\varphi_k\rangle\right)}\langle\varphi_k,M^{\alpha/2}_{-}\left(\mathbbm{1}_{[0,t)}u_{j,k}\right)\rangle\\
		&= \sum_{k=1}^d\frac{E_{\beta,\beta}\left(\frac{1}{2}\langle\varphi_k,\varphi_k\rangle\right)}{\beta E_\beta\left(\frac{1}{2}\langle\varphi_k,\varphi_k\rangle\right)}\int_{0}^{t} (M^{\alpha/2}_{+}\varphi_k)(s)\,u_{j,k}(s)\de s, \quad \varphi\in\mathcal{U}_0,
	\end{align*}
	 so that by Theorem~\ref{t:chrkd}, each component of the right-hand side of (\ref{e:SODE}) is holomorphic at zero in $\mathcal{S}_d(\R)_\C$ for all $t\in[0,T]$, and that
	\begin{align*}
		X_t = S_{\mu_{\beta}^{\otimes d}}^{-1}Y_t = V(t)V(0)^{-1}x_0 + \sigma V(t)\sum_{j=1}^{d}\langle\cdot,\eta_{t,j}\rangle\mathbf{e}_j.
	\end{align*}
	Finally, we show that the components of $(X_t)_{t\in[0,T]}$ satisfy the assumptions of Corollary~\ref{c:dvst}, that is, for some $p,q\in\N$, each component of the right-hand side of (\ref{e:ODE}) is uniformly bounded in $t\in [0,T]$ and $\varphi\in\mathcal{U}_{p,q}$. Now, using Estimate (\ref{e:CBound}), we can choose $p,q\in\N$ and a constant $K<\infty$ such that
	$$
		|C_{\beta,\alpha}(\varphi,t)|_{\text{euc}} \leq K,\quad t\geq 0,\ \varphi\in\mathcal{U}_{p,q}.
	$$
	Then the continuity of $A,V,V^{-1}$ on the compact interval $[0,T]$ yield a uniform bound of $|Y_t(\varphi)|_\text{euc}$, and hence, of $\left|\frac{\text{\rm d}}{\text{\rm d}t} Y_t(\varphi)\right|_\text{euc}$, in $t\in [0,T]$ and $\varphi\in\mathcal{U}_{p,q}$. Therefore, by Corollary~\ref{c:dvst}, we obtain the following result.
	\begin{thm}
		The process
		$$
			X_t = V(t)V(0)^{-1}x_0 + \sigma V(t)\sum_{j=1}^{d} \langle\cdot,\eta_{t,j}\rangle\mathbf{e}_j, \quad t\in [0,T],
		$$
		solves (\ref{e:DSDE}) as a system of equations in $(\mathcal{S}_d(\R))^{-1}_{\mu_{\beta}^{\otimes d}}$, where $V:[0,T]\to\R^{d\times d}$ is a fundamental matrix to the homogeneous system~(\ref{e:hODE}),
		$$
			\eta_{t,j} := \sum_{k=1}^{d}M^{\alpha/2}_{-}\left(\mathbbm{1}_{[0,t)}u_{j,k}\right)\mathbf{e}_k,
		$$
		and for $j,k=1,\dots,d$, $u_{j,k}:\R\to\R$ is the $(j,k)$-entry of $V^{-1}$, extended to zero outside $[0,T]$. Its $S_{\mu_{\beta}^{\otimes d}}$-transform is given by
		$$
			S_{\mu_{\beta}^{\otimes d}}X_t(\varphi) = V(t)V(0)^{-1}x_0 + \sigma V(t)\int_{0}^{t} V(s)^{-1} \sum_{k=1}^d\frac{E_{\beta,\beta}\left(\frac{1}{2}\langle\varphi_k,\varphi_k\rangle\right)}{\beta E_\beta\left(\frac{1}{2}\langle\varphi_k,\varphi_k\rangle\right)}\big(M^{\alpha/2}_{+}\varphi_k\big)(s)\,\mathbf{e}_k\de s,
		$$
		for $t\in [0,T]$ and $\varphi\in\mathcal{U}_0$, where $\mathcal{U}_0$ is a suitable neighborhood of zero in $\mathcal{S}_d(\R)_\C$.
	\end{thm}

\section*{Acknowledgments}	
The DAAD scholarship for K.~Orge within the Ph.D. program Mathematics in Industry and Commerce at TU Kaiserslautern is gratefully acknowledged.

\bibliographystyle{plain}

%\bibliography{masterref}

\end{document}